\documentclass[twocolumn]{autart}    % Enable this line and disable the 
\usepackage{ amsfonts, amssymb,textcomp}

\usepackage{amsmath}
\usepackage{enumitem}

\usepackage{graphicx}          % Include this line if your 
\newcommand{\al}{\alpha}
\newcommand{\be}{\beta}
\newcommand{\gm}{\gamma}

\newcommand{\dl}{\delta}

\newcommand{\s}{\sigma}

\newcommand{\Ree}{\operatorname{Re}}
\newcommand{\Imm}{\operatorname{Im\,}}

\newcommand{\T}{\operatorname{T\,}}

\newcommand{\id}{\operatorname{id}}
\newcommand{\supp}{\operatorname{supp\,}}
\newcommand{\spann}{\operatorname{span\,}}

\newtheorem{definition}{\hspace{\parindent}{\bf Definition}}
\newtheorem{lemma}{\hspace{\parindent}\bf Lemma}

\newtheorem{remark}{\hspace{\parindent}\bf Remark}

\newtheorem{assumption}{\hspace{\parindent}\bf Assumption}

\begin{document}

\begin{frontmatter}

\title{Stability analysis of coupled abstract systems with possibly unstable subsystem} % Title, preferably not more 

\thanks[footnoteinfo]{Some results introduced in this paper were partly presented at IEEE CDC 2024.\\ Corresponding author V.I. Slynko Tel. +380 67 316 49 47}

\author[Wue]{I.V. Atamas}\ead{ivan.atamas@uni-wuerzburg.de},   
\author[Wue]{D.Ye. Kulish} \ead{kulis6664@gmail.com},% Add the 
\author[Mech]{ A.A. Martynyuk}\ead{journalndst@gmail.com},  
\author[Mech]{V.I. Slynko}\ead{vitstab@ukr.net}  % (ead) as shown

\address[Wue]{Julius Maximilian University of W\"urzburg, Germany}  % Please supply                                              
\address[Mech]{Institute of Mechanics of NAS Ukraine, Ukraine}             % full addresses
       % here.

\begin{keyword}                           % Five to ten keywords,  
Lyapunov methods; 	Infinite-dimensional systems; 	Complex systems; Linear systems;  	Semigroup and operator theory           % chosen from the IFAC 
\end{keyword}                             % keyword list or with the 
                                          % help of the Automatica 
                                          % keyword wizard

\begin{abstract}                          
We consider a two-component coupled system of differential equations with operator coefficients. In contrast to the well-known small-gain approach, we assume the presence of the exponential stability property of only one subsystem. We introduce a new concept of  "dominance" of this subsystem, which allows us to find conditions on the exponential stability of a coupled system, and testing these on an example of infinite networks. The results are compared with a small-gain approach based on the Lyapunov vector function.
\end{abstract}

\end{frontmatter}

\section{Introduction}
% Математичне моделювання  процесів і явищ реального світу приводить до необхідності розглядати скінченновимірні \cite{Fieg} або нескінченно-вимірні системи диференціальних рівнянь у банаховому або гільбертовому просторах (див. \cite{ANM}   і бібліографію там). 
Mathematical modeling of real-world processes and phenomena leads to the need to consider finite-dimensional \cite{Fieg} or infinite-dimensional systems of differential equations in Banach or Hilbert spaces (see \cite{ANM} and bibliography there).

% Математичні моделі систем, які складаються зі скінченно-вимірних і нескінченновимірних компонент використовуються для опису процесів і явищ в фізиці, теорії керування, хімічній кінетиці та математичній біології. 
Mathematical models of systems consisting of finite-dimensional and infinite-dimensional components are used to describe processes and phenomena in physics, control theory, chemical kinetics, and mathematical biology.

% Математичні моделі хімічних і ядерних реакторів \cite{KaK19,Z2019}, напівпровідникових лазерів \cite{JR2001} описуються зв'язаними системами звичайних диференціальних рівнянь (ЗДР) і диференціальних рівнянь з частинними похідними (ДРЧП). Наприклад, в \cite{SMA2016} розглядається задача стабілізації крутильних і вісьових коливань, які виникають в роторних системах буріння, що моделюються зв'язаною системою ЗДР і ДРЧП. Зазначимо, що ефективні методи стабілізаціі таких систем розроблені в ряді робіт \cite{SCK2010}. Також в останні роки з розвитком різних типів мереж (соціальні мережі,транспортні і виробничі мережі) виникає проблема моделювання мереж великого розміру, які можна ідеалізувати як нескінченні мережі. Математичним образом таких мереж може бути нескінченна система ЗДР, яку можна подати в абстрактній формі як диференціальне рівняння або зв'язана система диференціальних рівнянь у банаховому просторі \cite{KM2022}.
Mathematical models of chemical and nuclear reactors \cite{KaK19,Z2019} as well as semiconductor lasers \cite{JR2001} are described by coupled systems of ordinary differential equations (ODEs) and partial differential equations (PDEs). For example, in \cite{SMA2016}, the problem of stabilizing torsional and axial vibrations that arise in rotary drilling systems modeled by a coupled system of ODEs and PDEs is considered. We note that effective methods of stabilization of such systems were developed in a number of works \cite{SCK2010}. Also, in recent years, with the development of various types of networks (social networks, transport and production networks), the problem of modeling networks of large size which can be idealized as infinite networks arose. The mathematical image of such networks can be an infinite system of ODEs which can be presented in an abstract form as a differential equation or a coupled system of differential equations in the Banach space \cite{KM2022}.

% Дослідження стійкості рівноважних станів в таких моделях є важливим етапом їх якісного аналізу. Лінеаризовані в околі стану рівноваги рівняння можна подати у вигляді двохкомпонентної абстрактної системи з операторними коефіцієнтами. Так само як і для скінченновимірних динамічних систем основним методом аналізу стійкості нескінченновимірних систем є метод функцій Ляпунова, заснований на скалярних \cite{Lyap}, векторних \cite{jpl}  або матрично-значних функціях \cite{aam}.
The study of stability of equilibrium states in such models is an important stage of their qualitative analysis. Equations linearized around the equilibrium state can be presented in the form of a two-component abstract system with operator coefficients. As well as for finite-dimensional dynamic systems, the main method of stability analysis of infinite-dimensional systems is the method of Lyapunov functions, based on scalar \cite{Lyap}, vector \cite{jpl} or matrix-valued functions \cite{aam}.

 % Метод функцій Ляпунова зручно використовувати тоді, коли  знаходження явних розв'язків системи є складною задачею. Важливо побудувати функцію Ляпунова яка б дозволила встановити з одного боку достатньо прості для перевірки умови стійкості. З іншого боку такі умови повинні бути не надто консервативними.
 The method of Lyapunov functions is convenient to use when finding explicit solutions of the system is a difficult task. On the one hand, it is important to construct a Lyapunov function that would allow establishing simple enough conditions for checking stability. On the other hand, such conditions should not be too conservative.

% Успіхи класичної і сучасної теорії зв'язних систем засновані на математичних методах, які надають можливість аналізувати більш складні структури виходячи з властивостей незалежних підсистем, кожна з яких володіє хорошими динамічними властивостями (наприклад, експоненціальної чи асимптотичної стійкості) \cite{AS2003}. Ця загальна парадигма є основою методів векторної функції Ляпунова (ВФЛ) \cite{belman, LMS1991} та small-gain теорем \cite{JTP1994,DRW2010}, на основі яких  отримано прості і досить ефективні достатні умови стійкості зв'язаних систем.
The success of the classical and modern theory of coupled systems are based on mathematical methods that make it possible to analyze more complex structures based on the properties of independent subsystems, each of which has good dynamic properties (for example, exponential or asymptotic stability) \cite{AS2003}. This general paradigm is the basis of the vector Lyapunov function (VLF) methods \cite{belman,LMS1991} and small-gain theorems \cite{JTP1994,DRW2010}, on the basis of which simple and quite effective sufficient conditions for the stability of coupled systems are obtained.

% В монографії \cite{ANM} наведено результати аналізу стійкості, нестійкості і обмеженості розв'язків диференціальних рівнянь в банаховому просторі на основі методів скалярних і векторних функцій Ляпунова. Як застосування розглядаються нейронні мережі з запізненням і модель багатоядерного атомного реактора.
The monograph \cite{ANM} gives the results of the analysis of stability, instability and boundedness of solutions of differential equations in the Banach space based on the methods of scalar and vector Lyapunov functions. Neural networks with a delay and a model of a multi-core nuclear reactor are considered as applications.

% При розгляді динамичних систем, динаміка яких описується різними типами рівнянь, а також зв'язаних великомасштабних  систем необходною умовою для застосування ВФЛ або small gain теорем є  асимптотична стійкість всіх незалежних підсистем, які утворюють складну систему. Недоліком цих методів є значна консервативність одержаних умов стійкості, оскільки припущення про експоненціальну чи асимптотичну стійкість незалежних підсистем не є необхідним для стійкості зв'язаної системи. Дослідження стійкості зв'язаних систем з нестійкими підсистемами потребують розробки принципово нових математичних підходів. Для дослідження стійкості складних динаміних систем з підсистемами, які володіють різними типами динамічної поведінки (тобто  є асимптотично стійкими,нейтрально стійкими і нестійкими) розроблено  метод матрично-значних функцій Ляпунова (МФЛ) (см. \cite{aam} і бібліографію там). Продуктивне використання цього  методу пов'язане з побудовою відповідної МФЛ.
When considering dynamic systems, the dynamics of which are described by various types of equations, as well as coupled large-scale systems \cite{djor83}, a necessary condition for the application of VLF or small gain theorems is the asymptotic stability of all independent subsystems that form a complex system. The disadvantage of these methods is the considerable conservatism of the obtained stability conditions, since the assumption of exponential or asymptotic stability of independent subsystems is not necessary for the stability of a coupled system. Studies of the stability of coupled systems with unstable subsystems require the development of fundamentally new mathematical approaches. To study the stability of complex dynamic systems with subsystems that have different types of dynamic behavior (i.e., asymptotically stable, neutrally stable, and unstable), the method of matrix-valued Lyapunov functions (MLF) was developed (see \cite{aam} and the bibliography there). The productive use of this method is coupled with the construction of the corresponding MLF.

% Для скінченно-вимірних динамічних систем з  дискретним часом в \cite{GL2015,NWZ2018a} запропоновано метод дослідження стійкості зв'язаних нелінійних систем який полягає в переході від вихідної дискретної системи до $M$ ітерованої дискретної системи. Для достатньо великих значень $M$ кількості ітерацій до $M$-ітерованої системи завжди застосовні класичні методи теорії стійкості зв'язанних систем, наприклад, small-gain теореми. Однак цей підхід поки що не отримав узагальнення для систем з неперервним часом. Інший метод дослідження стійкості лінійних зв'язаних періодичних систем базується на концепції матрично-значних функцій Ляпунова (МФЛ) розвивається в \cite{STA2022}. Нескінченно-вимірні зв'язані системи з нестійкими підсистемами розглядалися в \cite{DKS23, ADS2024}. В \cite{DKS23} розглянуто нелінійну зв'язану систему з вхідним сигналом, яка складається з нелійної системи ЗДР і рівняння з частинними похідними параболічного типу. Припускається, что система ЗДР не має властивості ISS, встановлено достатні умови ISS зв'язаної системи. В \cite{ADS2024} розглядається лінійна абстрактна система, яка складається з одновимірного ЗДР і абстрактного лінійного рівняння в гільбертовому простору. Припускається, що ЗДР може бути нестійким. Використовуючи теорію Крейна---Рутмана лінійних операторів в напівупорядкованих банахових просторах встановлено достатні умови експоненціальної стійкості.
In \cite{GL2015,NWZ2018a}, a method of studying the stability of coupled nonlinear systems for finite-dimensional dynamic systems with discrete time is considered. It consists of transitioning from the original discrete system to the $M$-iterated discrete system. For sufficiently large values of the number of iterations $M$, classical methods of stability theory of coupled systems are always applicable to an $M$-iterated system, for example, the small-gain theorem. However, this approach has not yet been generalized for continuous-time systems. Another method of studying the stability of linear coupled periodic systems is based on the concept of MLF developed in \cite{STA2022}. Infinite-dimensional coupled systems with unstable subsystems were considered in \cite{DKS23,ADS2024}. In \cite{DKS23}, a nonlinear coupled system with an input signal which consists of a nonlinear system of ODEs and an equation with partial derivatives of the parabolic type is considered. It is assumed that the ODE system does not have the ISS property, and sufficient conditions for the ISS coupled system are established. In \cite{ADS2024}, a linear abstract system which consists of a one-dimensional ODE and an abstract linear equation in Hilbert space is considered. It is assumed that the ODE may be unstable. Using the Krein---Rutman theory of linear operators in semi-ordered Banach spaces, sufficient conditions for exponential stability of the whole system are established.

% В цій статті пропонуєються нові методи побудови матричної функції Ляпунова (МФЛ) на основі яких встановлюються достатні умови стійкості зв'язаних двохкомпонентних лінійних систем з операторними коефіцієнтами. Припускається експоненціальна стікість лише однієї підсистеми за умови, що ця підсистема є домінуючою в сенсі запропонованої умови відділення спектрів.
This article proposes new methods for constructing the MLF based on which sufficient conditions for the stability of coupled two-component linear systems with operator coefficients are established. Exponential flow of only one subsystem is assumed, provided that this subsystem is dominant in the sense of the proposed spectrum separation condition.

% Стаття складається з семи розділів. У другому розділі  наводиться мотивцію дослідження і постановка проблеми. Третій розділ має технічний характер, у ньому обгрунтовано коректність постановки задачі. В четвертому розділі пропонується алгоритм побудови функції Ляпунова, на основі якого в розділі 5 доведено теореми про експоненціальну стійкість лінійних систем. В шостому розділі розглянуто приклади застосування доведених теорем: зчисленну систему диференціальних рівнянь і систему диференціальних рівнянь, яка складається зі звичайного диференціального рівняння і системи диференціальних рівнянь з частинними похідними. В останньому розділі обговорюються результати і перспективи подальших досліджень.
The article consists of eight sections. In the second section, the motivation of research and problem statement are given. The third section has a technical nature: it substantiates the well-posedness of the problem statement. In the fourth section, an algorithm for constructing the Lyapunov function is proposed. Based on this algorithm, theorems about the exponential stability of linear systems are proved in Section 5. In the sixth section, examples of application of proven theorems are considered: a countable system of differential equations and a system of differential equations which consists of an ordinary differential equation and a system of differential equations with partial derivatives. The last section discusses the results and prospects for further research.

{\it Notations.} If $M$ is a subset of some topological space, then $\overline M$ is its closure, $\id_M$ is the identical mapping $M\to M$. Let $H_i$, $i=1,2$ be complex Hilbert spaces with scalar product $\langle\cdot,\cdot\rangle_{H_i}$, which is linear in second component and semi-linear in the first, and norm $\|x\|_{H_i}=\sqrt{\langle x,x\rangle_{H_i}}$. The Banach space of linear bounded operators acting from $H_1$ to $H_2$ with the norm defined $\|A\|=\sup_{\|x\|_{H_1}=1}\|Ax\|_{H_2}$ is denoted by $L(H_1,H_2)$. The notation $L(H_1,H_1)\equiv L(H_1)$ is used. The conjugate space $H_1^{\prime}$ is identified with $H_1$ due to Ries's theorem. If $A\in L(H_1,H_2)$, then $A^{\dag}\in L(H_2,H_1)$ denotes a linear operator conjugate to $A$ , i.e., such that for all $x\in H_1$, $y\in H_2$ the identity $\langle y,Ax\rangle_{H_2}=\langle A^{\dag}y,x\rangle_{H_1}$ holds. For linear operator $A\,:\,\,\mathcal D(A)\to H_1$,  $\mathcal D(A)\subset H_1$ is its domain and hence a linear subspace, $\Ree A = 1/2(A+A^{\dag})$ is the real part of $A$, satisfying
\begin{equation*}
\Ree \langle x,Ax\rangle_{H_1}= \langle x,\Ree Ax\rangle_{H_1},
\end{equation*} for all $x \in \mathcal{D}(\Ree A)=\mathcal D(A) \cap \mathcal D(A^{\dag})$, and  $\rho(A)\subset\mathbb{C}$ is its resolvent set, i.e.,
\begin{equation*}
\rho(A):=\{\lambda\in\mathbb C\,:\, R_A(\lambda):=(\lambda\id_{H_1}-A)^{-1}\in L(H_1)\}.
\end{equation*} We call $A$ invertible, if $0$ belongs to the resolvent set, which is equivalent to operator $A^{-1}$ being bounded.
If $A\in L(H_1)$, $A^\dag=A$, then $A$ is called a self-adjoint operator, $\langle x,Ax\rangle_{H_1} \in \mathbb{R}, \, x\in H_i$ and we set
\begin{equation*}
\gathered
\lambda_{\min}(A):=\inf\limits_{\|x\|_{H_1}=1}\langle x,Ax\rangle_{H_1},\\
\lambda_{\max}(A):=\sup\limits_{\|x\|_{H_1}=1}\langle x,Ax\rangle_{H_1}.
\endgathered
\end{equation*}
In the case $\lambda_{\min}(A)>0$, the linear operator $A$ will be called positive definite.

We introduce Hilbert space $\ell^2:=\{(x_1,\dots,x_k,\dots)\,: x_k\in\mathbb C,\,k\in\mathbb N,\,\sum\limits_{k=1}^{\infty}|x_k|^2<+\infty\}$
with a scalar product $\langle x,y\rangle_{\ell^2}:=\sum\limits_{k=1}^{\infty}\overline{x_k}y_k$, $x$, $y\in\ell^2$.

The Hilbert space $L^2(0,\ell)$, $\ell>0$ consists of real-valued functions measurable on $(0,\ell)$, such that $\int_0^{\ell}f(x)^2 \,dx<+\infty$ with the scalar product $\langle f,g\rangle_{L^2(0,\ell)}:=\int_0^{\ell}f(x)g(x)\,dx$.

The Sobolev space $H^m(0,\ell)$ is defined (see \cite{Pa1983}) as the closure of the set $C^m[0,\ell]$ with respect to the scalar product
\begin{equation*}
\gathered
\langle f,g\rangle_{H^m(0,\ell)}=\sum\limits_{p=0}^m\langle f^{(p)},g^{(p)}\rangle_{L^2(0,\ell)}.
\endgathered
\end{equation*}
The Sobolev space $H^m_0(0,\ell)$ is defined as the closure of a set of finite functions
$C^m_0[0,\ell]=\{f\in C^m[0,\ell],\,\supp f\subset(0,\ell)\}$ with respect to the scalar product
$\langle f,g\rangle_{H^m(0,\ell)}$.\\
The tensor product of two vectors from some Hilbert space $H$ is defined by $a \otimes b : H \rightarrow \spann\{b\}, \,\, x \mapsto \langle a,x\rangle_{H}b$.

The derivative of a vector-valued function $ x : [0,T] \to H $ as well as the total derivative of the scalar function $V : H \to \mathbb{C} $, are defined in the standard way (see, for instance, \cite{krein-dal}, Chapter~1).\\

\section{Motivation.} Consider a linear coupled system of differential equations
\begin{equation}\label{3.1}
\gathered
\dot x(t)=Ax(t)+By(t),\quad x(0)=x_0\in H_1,\\
\dot y(t)=Cx(t)+Dy(t),\quad y(0)=y_0\in H_2,
\endgathered
\end{equation}
where $x\in H_1$, $y\in H_2$, $A\in L(H_1)$, $B\in L(H_2,H_1)$,
$C\in L(H_1,H_2)$, $D$ is a linear operator densely defined in $H_2$ (i.e., $\overline{\mathcal D(D)}=H_2$).

\begin{remark} Systems of linear differential equations written in the abstract form \eqref{3.1} include a wide class of linear infinite-dimensional systems. For example, coupled systems of ordinary differential equations and differential equations with partial derivatives of various types (parabolic, hyperbolic), coupled systems of integro-differential equations and equations with partial derivatives, etc.
\end{remark}

For an independent subsystem
\begin{equation}\label{3.2}
\gathered
\dot y(t)=Dy(t),\quad y(0)=y_0\in\mathcal D(D),
\endgathered
\end{equation}
we make the following assumption.
%The linear operator $D$ satisfies $0\in\rho(D)$ and
\begin{assumption} For some $\lambda_0 > 0$ operator $\lambda_0 \id_{H_2}-D$ is surjective and there is a (coercive) Lyapunov function $V_2(y)=\langle y,P_{2}y\rangle_{H_2}$, where $P_{2}\in L(H_2)$ is a linear self-adjoint operator such that for some positive constants $\widetilde{\alpha}_2$, $\widetilde{\beta}_2$, the inequalities
\begin{equation}\label{3.4}
\gathered
\widetilde{\alpha}_2\|y\|_{H_2}^2\le V_2(y)\le\widetilde{\beta}_2\|y\|_{H_2}^2
\endgathered
\end{equation}
hold and
the complete derivative of function $V_2(y)$ along the trajectory of the linear system \eqref{3.1} satisfies
\begin{equation}\label{3.4bis}
\gathered
\dot V_2(y)\Bigg|_{(2)}=\langle Dy,P_{2}y\rangle_{H_2}+\langle y,P_{2}Dy\rangle_{H_2}\\
\le-{\gm}_2\|y\|_{H_2}^2.
\endgathered
\end{equation}
for some ${\gm}_2 > 0$.
\end{assumption}
\begin{remark} It follows from Assumption 1 and the Lumer---Phillips Theorem (\cite{Pa1983}, p. 15, Corollary 4.3) that the linear operator $D$ generates a $C_0$-semigroup $(e^{tD})_{t\ge 0}$ in the Hilbert space $H_2$, and for positive constants $M_{_D} =\sqrt{\widetilde\beta_2/\widetilde\al_2}$, $\omega_{_D}=\gm_2/{2\widetilde\be_2}$ the inequality 
\begin{equation}\label{3.3}
\gathered
\|e^{tD}\|\le M_{_D}e^{-\omega_{_D} t},\quad t\ge 0.
\endgathered
\end{equation}
holds. In particular, $[0,\infty) \subset \rho(D)$.
\end{remark}
For the proof, we introduce a new scalar product $(x,y)_{H_2}:=\langle x,P_{2}y\rangle_{H_2}$ in the Hilbert space $H_2$.
Taking into account \eqref{3.4bis}, for all $y\in\mathcal D(D)$
\begin{equation}\label{3.3.1}
\gathered
\Ree (y,Dy)_{H_2}=
\frac{1}{2}\big(\langle Dy,P_{2}y\rangle_{H_2}+\langle y,P_{2}Dy\rangle_{H_2}\big)\leq0.
\endgathered
\end{equation}
We come to the conclusion that the linear operator $D$ is dissipative with respect to the scalar product $(.,.)_{H_2}$.
The norm $\|\cdot\|_2$ generated by this scalar product and the norm $\|\cdot\|_{H_2}$ are equivalent by virtue of \eqref{3.4}.
Considering Assumption 1 all the assumptions of the Lumer-Phillips Theorem are fulfilled, and we conclude that the linear operator $D$ generates a $C_0$-semigroup $(e^{tD})_{t\ge 0}$ in the Hilbert space $H_2$, which is contractive with respect to norm $\|\cdot\|_2$.\\ 
We derive the norm estimate \eqref{3.3} from the comparison theorem.\\ 
It follows $(-\omega_{_D},\infty) \subset \rho(D)$ from the Hille-Yosida Theorem (\cite{Pa1983}, p. 20, Theorem 5.3)
\\
We consider the problem of stability of the coupled system \eqref{3.1}, having previously recalled the definition of exponential stability.
\\Remark is proved.
\begin{definition} A linear system \eqref{3.1} is called exponentially stable if there are constants $M_0>0$, $\omega_0>0$ such that the inequality 
\begin{equation}\label{3.5}
\gathered
\|x(t)\|_{H_1}+\|y(t)\|_{H_2}\\ \le (\|x_0\|_{H_1}+\|y_0\|_{H_2})M_0e^{-\omega_0t},\quad t\ge 0
\endgathered
\end{equation}
holds.
\end{definition}

Various variants of the Lyapunov function method can be used to study the exponential stability of the \eqref{3.1} system.

The classical method of studying coupled systems, which arose in the 60s of the last century, is the method of Lyapunov vector functions (VLF). This method is based on the assumption of asymptotic or even exponential stability of independent subsystems. We assume the independent subsystem
\begin{equation}\label{3.6}
\gathered
\dot x(t)=Ax(t),\quad x(0)=x_0\in H_1
\endgathered
\end{equation}
in \eqref{3.1} to be exponentially stable. In this case, it is well known \cite{krein-dal} that for any linear self-adjoint positive definite operator $Q_1\in L(H_1)$, $Q_1^{\dag}=Q_1$ there exists a linear self-adjoint positive definite operator $ P_{1}\in L(H_1)$, $P_{1}^\dag=P_{1}$, which is a solution of the operator Lyapunov equation
\begin{equation}\label{3.7}
\gathered
A^{\dag}P_{1}+P_{1}A=-Q_1.
\endgathered
\end{equation}
For an independent subsystem \eqref{3.6}, one can construct a Lyapunov function $V_1(x)=\langle x,P_{1}x\rangle_{H_1}$, which has the following properties
\begin{equation}\label{3.8}
\gathered
\widetilde{\alpha}_1\|x\|_{H_1}^2\le V_1(x)\le\widetilde{\beta}_1\|x\|_{H_1}^2,\\
\dot V_1(x)=\langle x,(A^{\dag}P_{1}+P_{1}A)x\rangle_{H_1}\le-\gamma_1\|x\|_{H_1}^2,
\endgathered
\end{equation}
where $\widetilde{\alpha}_1=\lambda_{\min}(P_{1})$, $\widetilde{\beta}_1=\lambda_{\max}(P_{1})$ and
$ \gamma_1=\lambda_{\min}(Q_1)$.

% Для системи \eqref{3.1} на основі ВФЛ $V(x,y)=(V_1(x),V_2(y))^{\T}$,  стандартно встановлюються достатні умови експоненціальної стійкості лінійної системи \eqref{3.1} у вигляді нерівності
Based on the VLF $V(x,y)=(V_1(x),V_2(y))^{\T}$, sufficient conditions for the exponential stability of the linear system \eqref{3.1} are standardly established in the form of inequality
\begin{equation}\label{3.9}
\gathered
\|P_{1}B\|\|P_{2}C\|<\frac{ \gamma_1  \gamma_2}{4}.
\endgathered
\end{equation}
\begin{remark} A significant advantage of conditions \eqref{3.9} is their simplicity. However, they are conservative since the exponential stability of both independent subsystems is required.
\end{remark}

% Виникає питання як досліджувати експоненціальну стійкість у випадку, коли незалежна підсистема \eqref{3.6} не є експоненціально стійкою. Чи можна в цьому випадку отримати такі ж прості достатні умови стійкості зв'язаної системи \eqref{3.1} як і умови \eqref{3.9}. Для розв'язання цієї проблеми необхідно розвинути принципово нові підходи, відмінні від класичної парадигми ВФЛ. Для випадку $\dim H_1=1$ в \cite{ADS2024} використано теорію Крейна---Рутмана лінійних операторів у напівупорядкованих банахових просторах. Інший підхід, який буде розвинуто нижче, засновано на концепції матрично-значної функції Ляпунова (МФЛ) \cite{M1986}.
The question arises how to investigate exponential stability in the case if an independent subsystem \eqref{3.6} is not exponentially stable. Is it possible to obtain the same simple sufficient conditions for the stability of the coupled system \eqref{3.1} as \eqref{3.9}? To solve this problem, it is necessary to develop fundamentally new approaches, different from the classical paradigm of VLF. For the case $\dim H_1=1$ in \cite{ADS2024}, the Krein---Rutman theory of linear operators in semi-ordered Banach spaces can be used. Another approach, which will be developed below, is based on the concept of the matrix-valued Lyapunov function (MLF) \cite{mart-gut1979,dj1983,aam}.

We recall that the matrix-valued Lyapunov function is defined as a two-index system of functions
\begin{equation}\label{3.10}
\gathered
U(x,y)=
\begin{pmatrix}
v_{11}(x)&v_{12}(x,y)\\
v_{21}(x,y)& v_{22}(y)
\end{pmatrix},
\endgathered
\end{equation}
where  $v_{ii}\,:\,H_i\to\mathbb R$, $v_{ij}\,:\,H_1\times H_2\to\mathbb C$, $v_{ij}=\overline{v_{ji}}$, $i,j=1,2$. Based on MLF, it is possible to construct the scalar Lyapunov function $v(x,y,\eta)=\eta^{\T}U(x,y)\eta$, $\eta\in\mathbb R_{>0}^2$.

We note that if $v_{21}=v_{12}\equiv 0$, then MLF and VLF methods coincide. In this sense, the MLF method is a certain generalization of the VLF method. The question arises about the selection of MLF elements. In the case when the independent subsystem is exponentially stable, it is advisable to choose the corresponding element of the MLF as a Lyapunov function for the corresponding independent subsystem. For example, in our case $v_{22}(y)=V_2(y)=\langle y,P_{2}y\rangle_{H_2}$. It is necessary to select elements $v_{12}$ and $v_{11}$ of the MLF \eqref{3.10}. This paper is dedicated to solving these problems. To solve them, we introduce an important assumption regarding independent subsystems, which we call the domination condition.\\
\begin{assumption} (Domination condition) Let $D$ be invertible, then $D^{-k}:=(D^{-1})^k\in L(H_2)$, $k\in\mathbb N$. We assume, for some $r \in \mathbb{N}$ we can calculate the first $r-1$ positive powers of linear operator $A$, first $r$ powers of $D^{-1}$ and for the linear system \eqref{3.1} holds
\begin{equation}\label{3.11}
\|A^r\|\|D^{-r}\|\leq q_r <1.
\end{equation}If we can calculate all positive powers of operators $A$ and $D^{-1}$ and $$\inf\limits_{s\in\mathbb{N}}\|A^s\|\|D^{-s}\|<1$$ holds, we allow $r=\infty$ and set $q_r=0$. We call $r$ the computational limitation rate.
\end{assumption}
\begin{remark}
If Assumption 1 is already fulfilled, there is no need to check, whether the operator $D$ is invertible, since it follows from Remark 2.
\end{remark}
\if 0
\begin{remark}
Due to Hille-Yosida Theorem and estimate \eqref{3.3} the fulfillment of inequality
\begin{equation*}
\|A\|<\gm\sqrt{\frac{\al}{\be^3}},
\end{equation*} where $\al \leq \lambda_{\min}(P_{22}), \,\, \lambda_{\max}(P_{22})\leq \beta, \\\Ree\langle Dy, P_{22}y\rangle_{H_2} \leq -\gm \|y\|^2_{H_2}$ is sufficient for implementation of the {\it Domination condition}.
\end{remark}
\fi
\if 0
\begin{remark}
% Інтуїтивно умову відділення спектрів можна інтерпретувати так, що друга підсистема яка є експоненціально стійкою є домінуючою і має значний запас стійкості. Саме ця підсистема забезпечує стабілізацію зв'язаної системи при певних умовах на функції зв'язку між підсистемами.
Intuitively, the condition of separation of spectra can be interpreted so that the second subsystem, which is exponentially stable, is dominant and has a significant margin of stability. It is this subsystem that provides stabilization of the coupled system under certain conditions on the coupling function between subsystems.
\end{remark}
\fi

\section{Well-posedness of the problem} 
% В цьому розділі покажемо, що Припущення 1 гарантує коректність постановки задачі Коші \eqref{3.1}.Для цього використаємо операторний аналог формули Фробеніуса для обернення блочної матриці.
In this section, we will show that Assumption 1 guaranties the well-posedness of the formulation of the Cauchy problem \eqref{3.1}. For this, we will use the operator analog of the Frobenius formula for the inverse of a block matrix.

\begin{lemma}
Let
\begin{equation*}
\mathfrak{A}:=
\begin{pmatrix}
A&B\\
C&D
\end{pmatrix}
\end{equation*}
be a linear operator with domain $\mathcal D(\mathfrak{A})=H_1\times\mathcal D(D)\subset H_1\times H_2$ that acts according to the rule
\begin{equation*}
\mathfrak{A}\begin{pmatrix} x_1\\x_2\end{pmatrix}=
\begin{pmatrix} Ax_1+Bx_2\\Cx_1+Dx_2\end{pmatrix},\quad x_1\in H_1,\quad x_2\in\mathcal D(D).
\end{equation*}
% Позначимо $H:=A-BD^{-1}C$, і припустимо, що $D^{-1}\in L(H_2)$ і лінійний оператор $H\in L(H_1)$ оборотній, тоді лінійний оператор $\mathfrak{A}$ оборотній, тобто існує $\mathfrak{A}^{-1}\in L(H_1\times H_2)$ і справедлива операторна формула Фробеніуса%
We define scalar product on $H_1\times H_2$\,:
\begin{equation*}
\langle x,y\rangle_{H_1\times H_2}=\langle x_1,y_1\rangle_{H_1}+(x_2,y_2)_{H_2}
\end{equation*}
for $x=(x_1,x_2)\in H_1\times H_2$ and
$y=(y_1,y_2)\in H_1\times H_2$.\\
Let $K:=A-BD^{-1}C$ and suppose $D$ and $K$ both being invertible. Then the linear operator $\mathfrak{A}$ is invertible, and the Frobenius operator formula
\begin{equation*}
\mathfrak{A}^{-1}=
\begin{pmatrix}
K^{-1}&-K^{-1}BD^{-1}\\
-D^{-1}CK^{-1}&D^{-1}+D^{-1}CK^{-1}BD^{-1}
\end{pmatrix}
\end{equation*}
holds.
\end{lemma}

\begin{thm}
% Нехай для лінійної системи \eqref{3.1} виконуються умови Припущення 1, тоді задача Коші \eqref{3.1} при $t\in[0,\infty)$ генерує $C_0$-півгрупу на $H_1\times H_2$.
Let the conditions of Assumption 1 be satisfied for the linear system \eqref{3.1}. Then, the Cauchy problem \eqref{3.1} in $t\in[0,\infty)$ generates a $C_0$-semigroup in $H_1\times H_2$ and thus is well-posed.
\end{thm}

{\it Proof.} 
% Зауважимо, що не зменшуючи загальності, можна вважати, що скалярний добуток в просторі $H_2$ заданий формулою $(x,y)_{H_2}:=\langle x,P_{22}y\rangle_{H_2}$ (оскільки в силу \eqref{3.4} з Припущення 1 норма $\|x\|_2:=\sqrt{(x,x)_{H_2}}$ породжена скалярним добутком $(\cdot,\cdot)_{H_2}$ еквівалентна нормі $\|\cdot\|_{H_2}$). Нехай
We note that without loss of generality, we can assume that the scalar product in the space $H_2$ is given by the formula $(x,y)_{H_2}:=\langle x,P_{22}y\rangle_{H_2}$ by the equivalence of norms $\|\cdot\|_2$ and $\|\cdot\|_{H_2}$.\\ Let
\begin{equation*}
\mathfrak{A}:=
\begin{pmatrix}
A&B\\
C&D
\end{pmatrix}
\end{equation*}
where $A$, $B$, $C$ and $D$ are defined in \eqref{3.1}. 
Since $\overline{H_1\times\mathcal D(D)}=H_1\times\overline{\mathcal D(D)}=H_1\times H_2$, the linear operator $\mathfrak{A}$ is densely defined in $H_1\times H_2$. For $\omega\in\mathbb R$, we define the linear operator $\mathfrak{A}_{\omega}:=\mathfrak{A}-\omega\id_{H_1\times H_2}$, $\mathcal D(\mathfrak{A}_{\omega})=\mathcal D(\mathfrak{A})$. It suffices to prove that the linear operator  $\mathfrak{A}_{\omega}$ for some $\omega\in\mathbb R$ generates a $C_0$-semigroup in the space $H_1\times H_2$. For this, we will use the Lumer---Phillips theorem (Theorem 4.3 from \cite{Pa1983}, p. 14).
Taking into account \eqref{3.3.1}, for $z=(z_1,z_2)\in H_1\times\mathcal D(D)$ we obtain 
\begin{equation*}
\gathered
\Ree \langle z,\mathfrak{A}_{\omega}z\rangle_{H_1\times H_2}
\le(\|A\|-\omega)\|z_1\|_{H_1}^2\\+\big(2\|B\|+\beta_2\|P_{22}C\|\big)\|z_1\|_{H_1}\|z_2\|_{H_2}-\al_2\omega\|z_2\|_{H_2}^2.
\endgathered
\end{equation*} Continuing by using the Cauchy---Buniakovsky inequality we get
\begin{equation*}
\gathered
\Ree \langle z,\mathfrak{A}_{\omega}z\rangle_{H_1\times H_2} \le (\|A\|-\omega+\|B\|^2+\beta_2^2)\|z_1\|_{H_1}^2\\
+(\|P_{22}C\|^2-\alpha_2\omega)\|z_2\|_{H_2}^2
\endgathered
\end{equation*}
To guarantee the dissipativity of $\mathfrak{A}_{\omega}$, we take $\omega>\max\{\|A\|+\|B\|^2+\beta_2^2,\,\|P_{22}C\|^2/\alpha_2\}$.\\
Next, we show $\mathfrak{A}_{\omega}^{-1}\in L(H_1\times H_2)$. Remark 2 implies $\omega \in \rho(D)$, so, by Lemma 1 it is enough to check whether the operator
\begin{equation*}
\widehat K = \omega \id_{H_1}-A-BR_D(\omega)C
\end{equation*} is invertible, i.e. $\omega \in \rho\big(A+BR_D(\omega)C\big)$. Due to the Neumann Series, it suffices to show $\omega > \|A+BR_D(\omega)C\|$. By Hille-Yosida Theorem (\cite{Pa1983}, p. 20, Theorem 5.3) we get such estimate:
\begin{equation*}
\|A+BR_D(\omega)C\| \leq \|A\|+\|B\|\|C\|\frac{M_{_D}}{\omega+\omega_{_D}}.
\end{equation*} Thus, it is enough to choose $\omega > \|A\|+\|B\|\|C\|\frac{M_{_D}}{\omega_{_D}}$.
\if 
Next, we show that there exists a real number $\lambda_0>0$ such that $\lambda_0\in\rho(\mathfrak{A}_{\omega})$. Indeed,
$\lambda_0\id_{H_1\times H_2}-\mathfrak{A}_{\omega}=(\lambda_0+\omega)\id_{H_1\times H_2}-\mathfrak{A}$.
We denote $\widehat D:=(\lambda_0+\omega)\id_{H_2}-\mathfrak{A}$. 
As proved above, Assumption 1 guarantees that the linear operator $D$ generates a $C_0$-semigroup $(e^{tD})_{t\ge 0}$ in the space $H_2$. Therefore, for some $\beta\in\mathbb R$, we have $(\beta,+\infty)\subset\rho(D)$. 
For $\lambda_0>\beta-\omega$, There exists $\widehat D^{-1}\in L(H_2)$ and by the Hille---Yosida theorem (theorem 5.3 from \cite{Pa1983}, page 20) $\|\widehat D^{-1}\|\le\frac{M}{\lambda_0+\omega-\beta}$  for some $M>0$. According to the Frobenius operator formula, to conclude the proof, it is necessary to prove that the linear operator
\begin{equation*}
\gathered
\widehat H:=(\lambda_0+\omega)\id_{H_1}-A-B\widehat D^{-1}C\in L(H_1)
\endgathered
\end{equation*}
is invertible. If necessary, increase $\lambda_0$ so that
$\lambda_0>\|A\|-\omega$, then there exist $R_{A}(\lambda_0+\omega)\in L(H_1)$
and
\begin{equation*}
\gathered
\widehat H:=((\lambda_0+\omega)\id_{H_1}-A)(\id_{H_1}-R_{A}(\lambda_0+\omega)B\widehat D^{-1}C).
\endgathered
\end{equation*}
Since $A\in L(H_1)$, we have $\|R_{A}(\lambda_0+\omega)\|\le\frac{1}{\lambda_0+\omega-\|A\|}$. Therefore,
\begin{equation*}
\gathered
\|R_{A}(\lambda_0+\omega)B\widehat D^{-1}C\|\le
\frac{M\|B\|\|C\|}{(\lambda_0+\omega-\|A\|)(\lambda_0+\omega-\beta)}.
\endgathered
\end{equation*}
We increase $\lambda_0$ so that the inequality
\begin{equation*}
\gathered
(\lambda_0+\omega-\|A\|)(\lambda_0+\omega-\beta)>M\|B\|\|C\|,
\endgathered
\end{equation*}
 holds. Then, the linear operator  $(\id_{H_1}-R_{A}(\lambda_0+\omega)B\widehat D^{-1}C)$ is invertible, which proves the invertibility of
$\widehat H$. So, $\lambda_0\in\rho(\mathfrak{A}_{\omega})$. 
% На основі теореми Лумера---Філіпса приходимо до висновку що лінійний оператор $\mathfrak{A}$ генерує $C_0$-півгрупу в просторі $H_1\times H_2$. Теорему доведено.
\fi Based on the Lumer---Phillips theorem, we conclude that the linear operator $\mathfrak{A}$ generates a $C_0$-semigroup in the space $H_1\times H_2$. The theorem is proved.

\begin{remark} Let $(x(t),y(t))$ be a solution of the Cauchy problem \eqref{3.1} with initial conditions $(x(0),y(0))=(x_0,y_0)\in H_1\times \mathcal D(D)$. Then,
$(x(t),y(t))\in H_1\times \mathcal D(D)$ for all $t\ge 0$.
\end{remark}

\section{Construction of MLF} 
We choose the elements of the MLF $U(x,y)$ in the form
\begin{equation}\label{3.12}
\gathered
v_{11}(x)=\langle x,P_{11}x\rangle_{H_1},\quad v_{22}(y)=\langle y,P_{22}y\rangle_{H_2},\\
v_{12}(x,y)=\overline{v_{21}(x,y)}=\langle x, P_{12}y\rangle_{H_1},\\
\endgathered
\end{equation}
where $P_{ii}\in L(H_i)$, $i=1,2$ are linear self-adjoint operators and $P_{12}\in L(H_2,H_1)$ is adjointable.

% Обчислимо повні похідні елементів матрично-значної функції $U(x,y)$ вздовж розв'язків лінійної системи \eqref{3.1}, припускаючи, що $(x,y)\in H_1\times\mathcal D(D)$.
We calculate the complete derivatives of the elements of the matrix-valued function $U(x,y)$ along the solutions of the linear system \eqref{3.1} assuming that $(x,y)\in H_1\times\mathcal D(D)$\,:
\begin{equation}\label{3.13}
\gathered
\dot v_{11}(x)=2\Ree\langle x,P_{11}Ax\rangle_{H_1}
+2\Ree\langle x,P_{11}By\rangle_{H_1},\\
\dot v_{22}(y)=2\Ree\langle x,C^{\dag}P_{22}y\rangle_{H_1}+2\Ree\langle y,P_{22}Dy\rangle_{H_2},\\
\dot v_{12}(x,y)=\langle x,P_{12}Cx\rangle_{H_1}+\langle x,\big(A^{\dag}P_{12}+P_{12}D\big)y\rangle_{H_1}\\
+\langle y,B^{\dag}P_{12}y\rangle_{H_2}.
\endgathered
\end{equation}
% Не порушуючи загальності вектор $\eta$ виберемо у вигляді $\eta=(1)^{\T}$, $>0$, тоді для похідної скалярної функції $v(x,y,\eta)=\eta^{\T}U(x,y)\eta$ справедливе зображення
Without violating the generality, we choose the vector $\eta$ in the form $\eta=(1,\eta_0)^{\T}$, $\eta_0>0$. The corresponding scalar Lyapunov-function is given by
\begin{equation}\label{LF}
v(x,y,\eta_0)=v_{11}(x)+2 \eta_0 \Ree v_{12}(x,y) + \eta_0^2 v_{22}(y)
\end{equation}
Its derivative is
\begin{equation}\label{3.14}
\gathered
\dot v(x,y,\eta_0)=2\Ree \Big(
\langle x,\big(P_{11}A+\eta_0P_{12}C\big)x\rangle_{H_1}\\
+\langle x,\big(P_{11}B+\eta_0(A^{\dag}P_{12}+P_{12}D)+
\eta_0^2C^{\dag}P_{22}\big)y\rangle_{H_1}\\
+\eta_0\langle y,\big(B^{\dag}P_{12}+\eta_0P_{22}D\big)y\rangle_{H_2}\Big).
\endgathered
\end{equation}

We choose the linear operator $P_{12}$ as a solution of the operator equation
\begin{equation}\label{P_12_EQ}
\gathered
A^{\dag}P_{12}D^{-1}+P_{12}=-(\eta_0^{-1}P_{11}B
+\eta_0C^{\dag}P_{22})D^{-1}\\=:-PD^{-1}
\endgathered
\end{equation}
Due to the domination condition \eqref{3.11} $P_{12}$ can be found in the form
\begin{equation}\label{3.15}
\gathered
P_{12}=\sum\limits_{m=0}^{\infty}(-1)^{m+1}(A^{\dag})^{m}PD^{-(m+1)}\in L(H_2,H_1),
\endgathered
\end{equation}
% Права частина в \eqref{3.15} корректно визначена в силу умови відділення спектру \eqref{3.11} і справедлива оцінка  норми
with the norm estimate
\begin{equation}\label{3.16}
\gathered
\|P_{12}\|\le\frac{1}{1-q_r}\|P_{12}^{(r)}\|,
\endgathered
\end{equation}where 
\begin{equation*}
P_{12}^{(r)}=\sum\limits_{m=0}^{r-1}(-1)^{m+1}(A^{\dag})^{m}PD^{-(m+1)}.
\end{equation*}
% Формула \eqref{3.15} дозволяє спростити вираз для $\dot v(x,y,\eta)$, якщо врахувати, що
Substituting \eqref{3.15} in \eqref{LF} we find the scalar LF $v(x,y,\eta_0)$ as a linear mapping with respect to parameters $P_{11}$ and $\eta_0^2$, which is why, from now on, we can set $\eta_0=1$ without loss of generality as it leads to rescaling of the LF by factor $\eta_0^{-2}$.
The mid-term in \eqref{3.14} vanishes by virtue of \eqref{P_12_EQ} and we conclude
\if
The formula \eqref{3.15} allows to simplify the expression for $\dot v(x,y,\eta)$ if 
\begin{equation}\label{3.17}
\gathered
(A^{\dag}P_{12}y+P_{12}Dy)+P_{11}By+C^{\dag}P_{22}y=0
\endgathered
\end{equation}
is taken into account for all $y\in\mathcal D(D)$.
Taking into account \eqref{3.17}, the expression \eqref{3.14} takes the form
\fi
\begin{equation}\label{3.18}
\gathered
\dot v(x,y)=2\Ree \Big(
\langle x,\big(P_{11}A+P_{12}C\big)x\rangle_{H_1}\\
+\langle y,\big(B^{\dag}P_{12}+P_{22}D\big)y\rangle_{H_2}\Big).
\endgathered
\end{equation}
The following lemma establishes sandwich-estimates for the scalar function $v(x,y)$.

\begin{lemma} Let $\alpha_i\leq \lambda_{\min}(P_{ii})$, $\beta_i \geq \lambda_{\max}(P_{ii})$, $i=1,2$ and $q \geq \|P_{12}\|$ be arbitrary. For all $(x,y)\in H_1\times H_2$, the inequality
\begin{equation}\label{3.19}
\gathered
\lambda_{\min}(S_*)(\|x\|_{H_1}^2+\|y\|_{H_2}^2)\le v(x,y)\\\le\lambda_{\max}(S^*)(\|x\|_{H_1}^2+\|y\|_{H_2}^2),
\endgathered
\end{equation}
holds, where
\begin{equation}\label{3.20}
\gathered
S_*=
\begin{pmatrix}
\alpha_1&-q\\
-q&\alpha_2
\end{pmatrix},\,
S^*=
\begin{pmatrix}
\beta_1&q\\
q&\beta_2
\end{pmatrix}
\endgathered
\end{equation}
\end{lemma}

{\it Proof.} Indeed, taking into account \eqref{LF}, we find 
\begin{equation}\label{3.21}
\gathered
v(x,y)
\ge\lambda_{\min}(P_{11})\|x\|^2_{H_1}-2\|P_{12}\|\|x\|_{H_1}\|y\|_{H_2}\\
+\lambda_{\min}(P_{22})\|y\|_{H_2}^2
\ge\lambda_{\min}(S_*)(\|x\|_{H_1}^2+\|y\|_{H_2}^2),
\endgathered
\end{equation}
The upper bound for $v(x,y)$ is proved similarly. The lemma is proved.

We consider the problem of choosing the diagonal elements of the MLF. Since the independent subsystem \eqref{3.2} is exponentially stable by Remark 2, it is appropriate to choose $P_{22}=P_2$ (see Assumption 1). Since the independent subsystem \eqref{3.6} is not necessarily exponentially stable, a way of choosing the linear operator $P_{11}$ is proposed below.

Let $r\in\mathbb N\cup\{\infty\}$ be a computational limitation according to Assumption 2.
We consider the following linear operator inequality (LOI) for $P_{11}$\,:
\begin{equation}\label{3.22}
\gathered
G_r:=A^{\dag}P_{11}+P_{11}A\\
+\sum\limits_{m=0}^{r-1}(-1)^{m+1}\big((A^{\dag})^mP_{11}BD^{-(m+1)}C\\+C^{\dag}(D^{-(m+1)})^{\dag}B^{\dag}P_{11}A^m
\big)\\+\sum\limits_{m=0}^{r-1}(-1)^{m+1}((A^{\dag})^mC^{\dag}P_{22}D^{-(m+1)}C\\+C^{\dag}(D^{-(m+1)})^{\dag}P_{22}CA^m)<0.
\endgathered
\end{equation}

\if 0
{\it Second method.}  Let $Q_1$ be a self-adjoint positive definite linear operator, $r\in\mathbb Z_+\cup\{\infty\}$. 
For the linear operator $P_{11}$, we consider the linear operator equation 
\begin{equation}\label{3.27}
\gathered
A^{\dag}P_{11}+P_{11}A\\
+\sum\limits_{m=0}^{r}(-1)^{m+1}(C^{\dag}(D^{-(m+1)})^{\dag}B^{\dag}P_{11}A^m\\
+(A^{\dag})^mP_{11}BD^{-(m+1)}C\\
+(A^{\dag})^mC^{\dag}P_{22}D^{-(m+1)}C\\
+C^{\dag}(D^{-(m+1)})^{\dag}P_{22}CA^m)=-Q_1.
\endgathered
\end{equation}

% Формули (16) і операторне рівняння (23) (або (24)) вводяться так, щоб забезпечити якомога простіший для аналізу вигляд повної похідної $\dot v(x,y,\eta)$ функції $v(x,y,\eta)$ для лінійної системи (3.1).
Formulas \eqref{} and operator equation (23) (or (24)) are introduced in such a way as to provide the simplest form for analysis of the complete derivative $\dot v(x,y,\eta)$ of the function $v(x,y,\eta)$ for the linear system (3.1).
\fi

\section{Sufficient conditions for exponential stability}
% Лема 2 і формула \eqref{3.18} дозволяють сформулювати достатні умови експоненціальної стійкості лінійної системи \eqref{3.1}.
Lemma 2 and formula \eqref{3.18} make it possible to formulate sufficient conditions for the exponential stability of the linear system \eqref{3.1}.

{\bf Theorem 2.} 
% Припустимо, що для лінійної системи \eqref{3.1} виконуються умови Припущення 1 і умова відділення спектрів \eqref{3.11} і для  деякого $r\in\mathbb Z_+\cup\{\infty\}$ існують лінійний додатно-визначений самоспряжений оператор $P_{11}$, який задовольняє узагальнене операторне рівняння Ляпунова \eqref{3.22} і додатна стала $$ такі, що виконуються нерівності
Let Assumptions 1 and 2 be fulfilled for the linear system \eqref{3.1}, and for computational limitation $r\in\mathbb Z_+$, there exists a  positive definite linear operator $P_{11}$, which satisfies the LOI \eqref{3.22} such that the inequalities 
\begin{equation}\label{3.28}
\gathered
\|P_{12}^{(r)}\|<\sqrt{\alpha_1\alpha_2}(1-q_r),
\endgathered
\end{equation}
\begin{equation}\label{3.28bis}
\gathered
\frac{2q_r}{1-q_r}\|P_{12}^{(r)}\|\|C\|
<-\lambda_{\max}(G_r),
\endgathered
\end{equation}
\begin{equation}\label{3.28bisbis}
\gathered
\frac{2q_r}{1-q_r}\|P_{12}^{(r)}\|\|B\|<\gamma_2-\lambda_{\max}(F_r),
\endgathered
\end{equation}
are satisfied with the notation
\begin{equation*}
\gathered
P_{22}=P_2,\quad
\al_i\leq \lambda_{\min}(P_{ii}),\quad P=P_{11}B+C^{\dag}P_{22},\\
P_{12}^{(r)}=\sum\limits_{m=0}^{r-1}(-1)^{m+1}(A^{\dag})^{m}PD^{-(m+1)},\\
F_r= 2\Ree(B^{\dag}P_{12}^{(r)}).
\if 0 \sum\limits_{m=0}^r(-1)^{m+1}((B^{\dag}(A^{\dag})^mC^{\dag}P_{22}D^{-(m+1)}\\
+(D^{-(m+1)})^{\dag}P_{22}CA^mB)
+B^{\dag}(A^{\dag})^mP_{11}BD^{-(m+1)}\\
+(D^{-(m+1)})^{\dag}B^{\dag}P_{11}A^mB),\\ \fi
\endgathered
\end{equation*}
% Тоді лінійна система \eqref{3.1} експоненціально стійка.
Then, the linear system \eqref{3.1} is well-posed and exponentially stable.
\begin{remark}
% Теорема 2 зводить дослідження стійкості лінійної системи \eqref{3.1} до розв'язання операторного рівняння \eqref{3.22}, що є значно спрощує задачу. Наприклад, у випадку, коли $\dim H_1<+\infty$ рівняння \eqref{3.22} є лінійним матричним рівнянням, яке узагальнює рівняння Ляпунова і може бути розв'язане стандартними засобами лінійної алгебри, наприклад, на основі теорії прямих добутків матриць.
Theorem 2 reduces the study of the stability of linear system \eqref{3.1} to solving the LOI \eqref{3.22} which greatly simplifies the problem. For example, in the case $\dim H_1<+\infty$,  \eqref{3.22} is a linear matrix inequality which can be solved by standard means of linear algebra.
\end{remark}
\if 0
\begin{remark}
% Прямими обчисленнями на основі умови відділення спектрів \eqref{3.11} можна встановити оцінку
The following estimate 
\begin{equation*}
\gathered
\lambda_{\max}(F_r)\le
\frac{2(1-q^{r+1})\|D^{-1}\|\|B\|}{1-q}\|P_{11}B+C^{\dag}P_{22}\|.
\endgathered
\end{equation*}
can be established by direct calculations based on the spectrum separation condition \eqref{3.11}. 
% Ця нерівність дозволяє замінити останню нерівність системи \eqref{3.28} більш грубою нерівністю, перевага якої полягає в простішій перевірці
This inequality allows us to replace \eqref{3.28bisbis} by a rougher inequality, the advantage of which is easier verification
\begin{equation}\label{01}
\gathered
\frac{2\|D^{-1}\|\|B\|}{1-q}\|P_{11}B+C^{\dag}P_{22}\|<\gamma_2.
\endgathered
\end{equation}
\end{remark}
\fi

{\it Proof of Theorem 2.} Linear system \eqref{3.1} is well-posed by Theorem 1.\\
We consider the case $(x_0,y_0)\in H_1\times\mathcal D(D)$ first.
Then,  $(x(t),y(t))\in H_1\times\mathcal D(D)$ for $t\ge 0$.
By \eqref {3.18} and the elementary observation 
\begin{equation*}
\Ree\langle x,Sx\rangle \leq \lambda_{\max}(\Ree S)\|x\|^2,
\end{equation*} for arbitrary linear operator $S$ and vector $x$, we get the estimate
\begin{equation}\label{v_dot_est}
\gathered
\dot v(x,y)\leq 2\lambda_{\max}\Big(\Ree\big(P_{11}A+P_{12}C\big)\Big)\|x\|_{H_1}^2
\\+2\lambda_{\max}\big(\Ree(B^{\dag}P_{12})\big)\|y\|_{H_2}^2+2\Ree\langle P_{22}Dy,y\rangle
\endgathered
\end{equation}
We estimate the right-hand side in parts.\\
{1.} Notation in \eqref{3.22} is equivalent to 
\begin{equation*}
2\Ree (P_{11}A+P_{12}^{(r)}C)=G_r.
\end{equation*}
{2.} The rest of the $x$-term in \eqref{v_dot_est} we estimate using \eqref{3.11}\,:
\begin{equation*}
\gathered
\lambda_{\max}\Big(\Ree \big((P_{12}-P_{12}^{(r)})C\big)\Big)\\\leq \big\|P_{12}-P_{12}^{(r)}\big\|\|C\| \leq \frac{q_r}{1-q_r}\|P_{12}^{(r)}\|\|C\|.
\endgathered
\end{equation*}
{3.} Analogically, for $y$-term, 
\begin{equation*}
\gathered
\lambda_{\max}\big(2\Ree (B^{\dag}P_{12})\big)\leq \lambda_{\max}(F_r)\\+2\|B\|\big\|P_{12}-P_{12}^{(r)}\big\|\leq\lambda_{\max}(F_r)\\+\frac{2q_r}{1-q_r}\|P_{12}^{(r)}\|\|B\|,
\endgathered
\end{equation*}
and it follows by \eqref{3.4bis}, that
\begin{equation*}
\gathered
2\Ree\langle P_{22}Dy,y\rangle\leq -\gm_2\|y\|^2.
\endgathered
\end{equation*}
We conclude, that conditions \eqref{3.28bis}-\eqref{3.28bisbis} guarantee the existence of a positive constant $\mu>0$, for which
\begin{equation}\label{3.31}
\gathered
\dot v(x,y)\le-\mu(\|x\|_{H_1}^2+\|y\|_{H_2}^2)
\endgathered
\end{equation} holds. Coerciveness of function $v(x,y)$ is given by
condition \eqref{3.28}, which implies $\lambda_{\min}(S_*)>0$ when taking
\begin{equation*}
q = \frac{1}{1-q_r}\|P_{12}^{(r)}\| \geq \|P_{12}\|
\end{equation*} in Lemma 2, where the last estimate is given in \eqref{3.16}.\\
We continue by the classical exponential stability proof\,:\\
combining estimate \eqref{3.31} with the coerciveness of scalar LF $v(x,y)$, we establish 
\begin{equation}\label{3.32}
\gathered
\dot v(x,y)\le-\frac{\mu}{\lambda_{\max}( S^*)}v(x,y)
\endgathered
\end{equation}
Therefore, due to the comparison theorem and again \eqref{3.19}, we conclude 
\begin{equation}\label{3.33}
\gathered
\|x(t)\|_{H_1}^2+\|y(t)\|_{H_2}^2\le\frac{1}{\lambda_{\min}( S_*)}v(x(t),y(t))\\
\le
\frac{1}{\lambda_{\min}( S_*)}e^{-\frac{\mu}{\lambda_{\max}(S^*)}t}v(x_0,y_0)\\
\le\frac{\lambda_{\max}(S^*)}{\lambda_{\min}(S_*)}e^{-\frac{\mu}{\lambda_{\max}(S^*)}t}(\|x_0\|_{H_1}^2+\|y_0\|_{H_2}^2),\quad
t\ge 0,
\endgathered
\end{equation}
For {arbitrary} initial conditions, the validity of the inequality \eqref{3.33} follows from the density of a set $\overline{H_1\times\mathcal D(D)}$ in $H_1\times H_2$. The theorem is proved.

\section{One-dimensional case}
The motivation of this paragraph is that if $\dim H_1 = 1$, we can refuse the \textit{domination condition}.\\ 
Consider a coupled system of linear differential equations
\begin{equation}\label{3.46*}
\gathered
\dot x=ax+\langle b,y\rangle_{H},\quad x(0)=x_0\in\mathbb C,\\
\dot y=cx+Dy,\quad y(0)=y_{0}\in H,
\endgathered
\end{equation}
where $a,x\in\mathbb C$, $b,c\in H$, $\overline{\mathcal D(D)}=H$, $\dim H \geq 2$. %,$\, \forall y\in \mathcal D(D)$.% 
\begin{assumption} The linear operator $D$ satisfies $-\bar a \in\rho(D)$, for some $\lambda_0 \geq 0$ operator $\lambda_0 \id_H-D$ is surjective and there is a positive constant $\gm_0$ satisfying
\begin{equation}\label{3.46**}
\gathered
\Ree\langle y,Dy\rangle\le-\gamma_0\|y\|^2, \quad y \in \mathcal{D}(D).
\endgathered
\end{equation}
\end{assumption}
We denote 
\begin{equation*}
\gathered
R_{D}(-\overline a)=-(\overline a\id+D)^{-1}=-D^{-1}(\id+\overline aD^{-1})^{-1},\\
\rho_b=\|R_{D}(-\overline a)^{\dag}b\|, \quad \rho_c=\|R_{D}(-\overline a)^{\dag}c\|,\\
\rho_{cc}=\langle R_{D}(-\overline a)c,c\rangle,\quad
\rho_{bb}=\langle R_{D}(-\overline a)b,b\rangle,\\
\rho_{bc}=\langle R_{D}(-\overline a)c,b\rangle,\quad
{\rho}_{abc}:=\Ree(a+\rho_{bc}),
\endgathered
\end{equation*}and
\begin{equation*}
\gathered
\mu_{bc}:=\langle R_{D}(-\overline a)^{\dag}b,R_{D}(-\overline a)^{\dag}c\rangle,\\
\chi_b:=\frac{\rho_b^2}{1-2\Ree\mu_{bc}},\quad
\chi_c:=\frac{\rho_c^2}{1-2\Ree\mu_{bc}},
\endgathered
\end{equation*}
\begin{equation*}
\gathered
\theta_{\pm}:=\frac{1\pm\sqrt{1-4\chi_b\chi_c}}{2\chi_b}
\endgathered
\end{equation*}
\begin{equation*}
\gathered
%\omega_{bc}=\rho_c\|b\|_H+\rho_{bc},\\
%\omega_{bb}=\rho_b\|b\|_H+\rho_{bb}
\omega_{bc}:=\Ree\rho_{bc}+\sqrt{\rho_c^2\|b\|_H^2-(\Imm\rho_{bc})^2},\\
\omega_{bb}:=\Ree\rho_{bb}+\sqrt{\rho_b^2\|b\|_H^2-(\Imm\rho_{bb})^2},\\
\endgathered
\end{equation*}

{\bf Theorem 3.} Let linear system \eqref{3.46*} satisfy Assumption 3 and inequalities
\begin{equation*}
\gathered
2\Ree\mu_{bc}<1,\quad
4\chi_b\chi_c < 1, \quad \theta_-\omega_{bb}<2\gm_0-\omega_{bc}.
\endgathered
\end{equation*}
Fulfillment of one of the following systems of inequalities

\begin{enumerate}[label=\roman*)]%\label{T3_Cond}
\item       ${\rho}_{abc}<0,\quad \Ree\rho_{cc}<\theta_+|{\rho}_{abc}|,\\ \\
\omega_{bb}\Ree\rho_{cc}
<(2\gm_0-\omega_{bc})|{\rho}_{abc}|$,\\
\item  ${\rho}_{abc}=0,\quad \Ree\rho_{cc}<0$\\
\item ${\rho}_{abc}>0,\quad \Ree\rho_{cc}<\theta_-{\rho}_{abc}$
\end{enumerate}
\if 0
The existence of some constant  $\eta_0>0$ such that
\begin{equation}\label{SC1}
\gathered
\frac{1+\Ree\rho_{cc}}{\Ree(a+\rho_{bc})}<0,
\endgathered
\end{equation}
\begin{equation}\label{SC2}
\gathered
\pm\Big(1+\frac{\Ree\rho_{cc}\theta_{\pm}}{\Ree(a+\rho_{bc})}\Big)\pm\frac{\theta_{\pm}}{\Ree(a+\rho_{bc})}<0,
\endgathered
\end{equation}
\begin{equation}\label{SC3}
\gathered
\Big(2\gamma_0-\omega_{bc}+\frac{\omega_{bb}\Ree\rho_{cc}}{\Ree(a+\rho_{bc})}\Big)
+\frac{\omega_{bb}}{\Ree(a+\rho_{bc})}>0.
\endgathered
\end{equation}
\fi
guarantees the well-posedness exponential stability of the linear system \eqref{3.46*}.

\section{Examples and comparison of results}

\subsection{Exponential $\ell^2$-stability of a countable system}
% Розглянемо  зчисленну систему лінійних диференціальних рівнянь
Consider a countable system of linear differential equations
\begin{equation}\label{3.46}
\gathered
\dot x=ax+\sum\limits_{k=1}^{\infty}\overline b_ky_k,\quad x(0)=x_0\in\mathbb C,\\
\dot y_n=c_nx-\lambda_n y_n,\quad y_n(0)=y_{n0}\in\mathbb C,\quad n\in\mathbb N.
\endgathered
\end{equation}
Suppose that $y_0=(y_{n0})_n, b:=(b_n)_n, c:=(c_n)_n$ are vectors in $\ell^2$ and $(\lambda_n)_n\subset\mathbb C$ satisfies
$\inf\limits_{n}(\Ree\lambda_n)>0$ and $\inf\limits_{n}(|\lambda_n-\overline a|)>0$.

Let $D$ be a linear operator with the domain
\begin{equation*}
\mathcal{D}(D)=\{y\in\ell^2\,\,:\,\,\sum\limits_{k=1}^{\infty}|\lambda_k|^2|y_k|^2<+\infty\}
\end{equation*}
which acts according to the rule
\begin{equation*}
Dy=-(\lambda_1y_1,\dots,\lambda_ny_n,\dots).
\end{equation*}
% Лінійну зчисленну систему \eqref{3.46} подамо в абстрактному вигляді як зв'язану систему на добутку просторів
Linear operator $D$ is densely defined since its domain contains a set of finite sequences, which is dense in $\ell^2$.\\
With such operator $D$, the countable linear system \eqref{3.46} can be represented in the form of a coupled system \eqref{3.46*}.\\
% Зауважимо, що експоненціальна стійкість системи \eqref{3.47} в сенсі означення 1 еквівалентно експоненціальній $\ell^2$-стійкості зчисленної системи \eqref{3.46}.
We note that the exponential stability of the system \eqref{3.46} in the sense of Definition 1 is equivalent to the exponential $\ell^2$-stability of the countable system \eqref{3.46}.

We verify the conditions of Assumption 3. For an arbitrary $z \in \ell^2$ and $\lambda \in \mathbb{C}$, satisfying $\inf\limits_{n}(|\lambda_n+\lambda|)>0$, it is easy to see,
\begin{equation*}
R_D(\lambda)z=\big(\frac{z_1}{\lambda+\lambda_1},\dots,\frac{z_n}{\lambda+\lambda_n},\dots\big) \in \ell^2.
\end{equation*} It follows $-\overline{a},1 \in \rho(D)$. For all $y\in\mathcal D(D)$,
\begin{equation}\label{3.48}
\gathered
\Ree\langle y,Dy\rangle_{\ell^2}=-\sum\limits_{k=1}^{\infty} \Ree \lambda_k|y_k|^2 \le-\inf\limits_{k}(\Ree\lambda_k)\|y\|_{\ell^2}^2.
\endgathered
\end{equation}
Therefore, the conditions of Assumption 3 are fulfilled with $\gamma_0=\inf\limits_n(\Ree\lambda_n)>0$ and system \eqref{3.46} is well-posed by the virtue of Theorem 1.

Based on the VLF method, from \eqref{3.9}, we obtain the conditions of exponential stability of the coupled system
\begin{equation}\label{3.49}
\gathered
\Ree a<0,\quad \|b\|_{\ell^2}\|c\|_{\ell^2}<-\gm_0\Ree a.
\endgathered
\end{equation}

We define the constants from Theorem 3\,:
\begin{equation}\label{3.51*}
\gathered
\rho_b=\bigg(\sum\limits_{k=1}^{\infty}\frac{|b_k|^2}{|\lambda_k-\overline a|^2}\bigg)^{1/2},\quad
\rho_c=\bigg(\sum\limits_{k=1}^{\infty}\frac{|c_k|^2}{|\lambda_k-\overline a|^2}\bigg)^{1/2},\\
\rho_{cc}=\sum\limits_{k=1}^{\infty}\frac{|c_k|^2}{\overline \lambda_k- a},\quad
\rho_{bb}=\sum\limits_{k=1}^{\infty}\frac{|b_k|^2}{\overline \lambda_k- a},\\
\quad\rho_{bc}=\sum\limits_{k=1}^{\infty}\frac{
 b_k\overline c_k}{\lambda_k-\overline a},\quad 
\mu_{bc}=\sum\limits_{k=1}^{\infty}\frac{
\overline b_kc_k}{|\lambda_k-\overline a|^2}
\endgathered
\end{equation} to apply Theorem 3.

We consider a numerical example choosing $a=0.2$, $\lambda_k=k^2$, $b:=\tau\widehat b$, $c:=\tau\widehat c$, coefficient of amplification of couple between subsystems $\tau>0$,
\begin{equation*}
\gathered
\widehat{b}:=\big(1,\dots,\frac{1}{k},\dots\big)\in\ell^2,\quad
\widehat{c}:=-\big(1,\dots,\frac{1}{k^{3/2}},\dots\big)\in\ell^2.
\endgathered
\end{equation*}
For these parameters, $\|\widehat{b}\|_{\ell^2}=1.2825$,
$\rho_{\widehat{b}}=1.2576$, 
$\rho_{\widehat{c}}=1.2537$,
$\rho_{\widehat{b}\widehat{b}}=1.336$,
$\rho_{\widehat{c}\widehat{c}}=1.2887$,
$\rho_{\widehat{b}\widehat{c}}=-1.3072$.

% Твердження 1 (Теорема 2) гарантує експоненціальну $\ell^2$-стійкість лінійної системи \eqref{3.46}  для значень коефіцієнту підсилення $\tau\in(0.5728,0.81)$. Оскільки $a>0$, перша підсистема є нестійкою, тому метод ВФЛ не застосовний в цьому випадку.
Theorem 3 guarantees exponential $\ell^2$-stability of the linear system \eqref{3.46} for the values of the gain coefficient $\tau\in(0.4893,0.6312)$. Since $a>0$, the first subsystem is unstable, so the VLF method is not applicable in this case.

For the parameter value $a=-0.25$ and the same $\lambda_k$, $b$, $c$, we find
$\rho_{\widehat{b}}=0.8428$, 
$\rho_{\widehat{c}}=0.8379$,
$\rho_{\widehat{b}\widehat{b}}=0.9123$,
$\rho_{\widehat{c}\widehat{c}}=0.8687$,
$\rho_{\widehat{b}\widehat{c}}=-0.8857$.
% Твердження 1 (Теорема 2) гарантує експоненціальну $\ell^2$-стійкість лінійної системи \eqref{3.46}  для значень коефіцієнту підсилення $\tau\in(0,0.9572)$. Оскільки $a<0$, перша підсистема є  експоненціально стійкою, тому можна використати метод ВФЛ і умови експоненціальної стійкості \eqref{3.9}. Ці умови приводять до значень коефіцієнту підсилення з інтервалу $(0,0.4216)$.
Theorem 3 guarantees exponential $\ell^2$-stability of the linear system \eqref{3.46} for the values of the gain coefficient $\tau\in(0,1.0507)$. Since $a<0$, the first subsystem is exponentially stable, so it is possible to use the VLF method and the conditions of exponential stability \eqref{3.9}. These conditions lead to $\tau \in (0,0.37714)$.

\subsection{A coupled ODE+PDE system of parabolic type}
Consider a coupled system of the form
\begin{equation}\label{3.65}
\gathered
\dot x(t)=ax(t)+\int\limits_{0}^{\ell} b(z) y(z,t)\,dz,\\
\frac{\partial y(z,t)}{\partial t}= c(z)x(t)+\frac{\partial }{\partial z}\Big(d(z)\frac{\partial y(z,t)}{\partial z}\Big),
\endgathered
\end{equation}
with initial and boundary conditions
$x(0)=x_0\in\mathbb R$, $y(z,0)=y_0(z)\in H^2(0,\ell)\cap H_0^1(0,\ell)$,
$ y(0,t)=y(\ell,t)=0$, $(z,t)\in[0,\ell]\times \mathbb R_+$. Here $\ell>0$, $a\in\mathbb R$, $b$, $c\in L^2(0,\ell)$,
$d\in C^1[0,\ell]$, $d(z)\ge d_0>0$, $d_0=\min_{z\in[0,\ell]}d(z)$. This system takes the form \eqref{3.46*} when defining the linear operator $D\,:\,\,\mathcal{D}(D)\to L^2(0,\ell)$,
\begin{equation*}\gathered
\mathcal{D}(D)=H^2(0,\ell)\cap H_0^1(0,\ell),\quad 
(Dy)(z)=\frac{d}{dz}\Big(d(z)\frac{dy}{dz}(z)\Big).
\endgathered
\end{equation*}
This operator is densely defined in $L^2(0,\ell)$, and also, due to the Friedrichs inequality, we have
\begin{equation*}
\gathered
\langle y,Dy\rangle=\int\limits_0^{\ell}y(z)Dy(z)\,dz\\
=-\int\limits_0^{\ell}d(z)y'(z)^2\,dz
\le -\frac{\pi^2d_0}{\ell^2}\|y\|^2.
\endgathered
\end{equation*}
Based on the VLF method, we obtain the following conditions of exponential stability\,:
\begin{equation}\label{3.499}
\gathered
a<0,\quad \|b\|\cdot\|c\|<-\frac{\pi^2d_0a}{\ell^2}.
\endgathered
\end{equation}
We will not apply Theorem 3 to the linear system \eqref{3.65}, since it is hard to find the resolvent $R_D(-a)$ for this example -- we check the conditions of Assumption 1 first. Let $P_{2}=\id_{L^2(0,\ell)}$, $\gamma_2={2\pi^2d_0}/{\ell^2}$. The linear differential equation $Dy=f$ is solvable for all $f\in L^2(0,\ell)$\,:
\begin{equation*}
\gathered
y(z)=(D^{-1}f)(z)=\int\limits_0^{\ell}\Gamma(z,\zeta)f(\zeta)\,d\zeta,
\endgathered
\end{equation*}
where $\Gamma$ is the Green's function
\if 0
\begin{equation*}
\gathered
\Gamma(z,\zeta)=\begin{cases}
-\frac{\int\limits_0^{z}\frac{d\xi}{d(\xi)}\int\limits_{\zeta}^{\ell}\frac{d\xi}{d(\xi)}}{\int\limits_0^{\ell}\frac{d\xi}{d(\xi)}},\quad 0\le z\le\zeta\le\ell\\
-\frac{\int\limits_0^{\zeta}\frac{d\xi}{d(\xi)}\int\limits_{z}^{\ell}\frac{d\xi}{d(\xi)}}{\int\limits_0^{\ell}\frac{d\xi}{d(\xi)}},\quad 0\le \zeta\le z\le\ell
\end{cases}
\endgathered
\end{equation*}
\fi
\begin{equation*}
\gathered
\Gamma(z,\zeta)=-\frac{\int\limits_0^{\min\{z,\zeta\}}\frac{d\xi}{d(\xi)}\int\limits_{\max\{z,\zeta\}}^{\ell}\frac{d\xi}{d(\xi)}}{\int\limits_0^{\ell}\frac{d\xi}{d(\xi)}}, \quad z,\zeta \in [0,\ell], 
\endgathered
\end{equation*}
satisfying additionally $\Gamma(z,\zeta)=\Gamma(\zeta,z)$, which verifies that operator $D^{-1}$ is self-adjoint.
%Thus,  $0\in\rho(D)$ and Assumption 1 holds.

% Визначимо ітеровані ядра
We define iterated kernels
\begin{equation*}
\gathered
\Gamma_1(z,\zeta):=\Gamma(z,\zeta),\\
\Gamma_{k+1}(z,\zeta):=\int\limits_0^{\ell}\Gamma_k(z,s)\Gamma(s,\zeta)\,ds,\quad k\in\mathbb N,
\endgathered
\end{equation*} which satisfy 
\begin{equation*}
\gathered
\big(D^{-k}f\big)(z)=\int\limits_0^{\ell}\Gamma_k(z,\zeta)f(\zeta)\,d\zeta
\endgathered
\end{equation*} due to the Fubini Theorem. Thus, $D^{-k}$ is an integral operator, which implies
\begin{equation*}
\gathered
\|D^{-k}\|^2=\int\limits_{0}^{\ell}\,\int\limits_{0}^{\ell}\Gamma_{k}^2(z,\zeta)\,d\zeta\,dz.
\endgathered
\end{equation*}
In particular, $0\in\rho(D)$ and so, fulfillment of Assumption 1 follows from the openness of the resolvent set. Let Assumption 2 be fulfilled for finite computational limitation $r$.\\
We denote 
\begin{equation*}
\gathered
\Psi_r(z,\zeta)=\sum\limits_{k=0}^{r-1}(-1)^{k+1}a^k\Gamma_{k+1}(z,\zeta),\\
\rho_b^{(r)}=\frac{1}{1-q_r}\Big(\int\limits_0^{\ell}\Big(\int\limits_0^{\ell}\Psi_r(z,\zeta)b(\zeta)\,d\zeta\Big)^2\,dz\Big)^{1/2}\\
\rho_c^{(r)}=\frac{1}{1-q_r}\Big(\int\limits_0^{\ell}\Big(\int\limits_0^{\ell}\Psi_r(z,\zeta)c(\zeta)\,d\zeta\Big)^2\,dz\Big)^{1/2}\\
\rho_{cc}^{(r)}=q_r\|c\|\rho_c^{(r)}+\int\limits_0^{\ell}\int\limits_0^{\ell}\Psi_{r}(z,\zeta)c(z)c(\zeta)\,dz\,d\zeta\\
\rho_{bb}^{(r)}=q_r\|b\|\rho_b^{(r)}+\int\limits_0^{\ell}\int\limits_0^{\ell}\Psi_{r}(z,\zeta)b(z)b(\zeta)\,dz\,d\zeta,\\
\rho_{bc}^{(r)}=q_r\|b\|\rho_c^{(r)}+\int\limits_0^{\ell}\int\limits_0^{\ell}\Psi_{r}(z,\zeta)c(z)b(\zeta)\,dz\,d\zeta,\\
{\rho}_{abc}^{(r)}:=a+\rho_{bc}^{(r)},
\endgathered
\end{equation*}and
\begin{equation*}
\gathered
\mu_{bc}^{(r)}:=q_r(2-q_r) \rho_b^{(r)} \rho_c^{(r)}\\+\int\limits_0^{\ell}\int\limits_0^{\ell}\int\limits_0^{\ell}\Psi_r(z,\zeta)
\Psi_r(z,\xi)b(\zeta)c(\xi)\,dz\,d\xi\,d\zeta,\\
\chi_b^{(r)}:=\frac{(\rho_b^{(r)})^2}{1-2\mu_{bc}^{(r)}},\quad
\chi_c^{(r)}:=\frac{(\rho_c^{(r)})^2}{1-2\mu_{bc}^{(r)}},
\endgathered
\end{equation*}
\begin{equation*}
\gathered
\theta_{\pm}^{(r)}:=\frac{1\pm\sqrt{1-4\chi_b^{(r)}\chi_c^{(r)}}}{2\chi_b^{(r)}}
\endgathered
\end{equation*}
\begin{equation*}
\gathered
\omega_{bc}^{(r)}={\rho}_{bc}^{(r)}+\|b\|{\rho}_{c}^{(r)},\quad
\omega_{bb}^{(r)}={\rho}_{bb}^{(r)}+\|b\|{\rho}_{b}^{(r)}.
\endgathered
\end{equation*}
{\bf Proposition 1.} Let linear system \eqref{3.46*} satisfy Assumption 2 and inequalities
\begin{equation*}
\gathered
2\mu_{bc}^{(r)}<1,\quad
4\chi_b^{(r)}\chi_c^{(r)} < 1, \quad \theta_-^{(r)}\omega_{bb}^{(r)}<\frac{2\pi^2d_0}{\ell^2}-\omega_{bc}^{(r)}.
\endgathered
\end{equation*}
Fulfillment of one of the following systems of inequalities
\begin{enumerate}[label=\roman*)]%\label{T3_Cond}
\item ${\rho}_{abc}^{(r)}<0,\quad \rho_{cc}^{(r)}<\theta_+^{(r)}|{\rho}_{abc}^{(r)}|,\\ \\ \omega_{bb}^{(r)}\rho_{cc}^{(r)}<(\frac{2\pi^2d_0}{\ell^2}-\omega_{bc}^{(r)})|{\rho}_{abc}^{(r)}|$,\\
\item ${\rho}_{abc}^{(r)}=0,\quad \rho_{cc}^{(r)}<0$,\\
\item ${\rho}_{abc}^{(r)}>0,\quad \rho_{cc}^{(r)}<\theta_-^{(r)}{\rho}_{abc}^{(r)}$.
\end{enumerate}
Then, the linear system \eqref{3.65} is exponentially stable.
{\it Proof} is given in Appendix B.

We consider a numerical example of a linear system \eqref{3.65} with parameters $\ell=1$, $d(z)=\frac{1}{1+z^2}$, $b(z)=\tau z$, $c(z)=-\tau z^3$ with a connection amplification factor $\tau>0$ .
The Green's function $\Gamma(z,\xi)$ is given by
\begin{equation*}
\gathered
\Gamma(z,\xi)=
\begin{cases}
\frac{1}{12}z(z^2+3)(\xi^3+3\xi-4),\quad 0\le z\le\xi\le 1,\\
\frac{1}{12}\xi(\xi^2+3)(z^3+3z-4),\quad 0\le\xi\le z\le 1
\end{cases}
\endgathered
\end{equation*}
We first consider $a=0.2$. Proposition 1, applied with the calculation limitation $r=3$ guarantees exponential stability for the values of the gain coefficient $\tau\in[3.75,16.068]$. Since $a>0$, the first subsystem is unstable and the VLF method is not applicable in this case.

In case of $a=-0.2$, we get the confidence intervall  $\tau \in [0,2.127)$ via VLF method solving the inequality \eqref{3.499} with respect to the gain coefficient $\tau$, while Proposition 1, applied with the calculation limitation $r=4$ guarantees exponential stability for the values of the gain coefficient $\tau\in[0,17.265]$.
% Твердження 1 (Теорема 2) гарантує експоненціальну $\ell^2$-стійкість лінійної системи \eqref{3.46}  для значень коефіцієнту підсилення $\tau\in(0,0.9572)$. Оскільки $a<0$, перша підсистема є  експоненціально стійкою, тому можна використати метод ВФЛ і умови експоненціальної стійкості \eqref{3.9}. Ці умови приводять до значень коефіцієнту підсилення з інтервалу $(0,0.4216)$.

\section{Conclusion} 
% Основним результатом доповіді є явна формула (16) для визначення поза-діагональних елементів МФЛ і операторні рівняння (23) (або (24)) які визначають діагональний елемент, що відповідає, можливо,  нестійкій підсистемі. Ці результати розвивають результати статті [24] для більш широкого класу нескінченно-вимірних динамічних систем.
The main result of the report is the explicit formula (16) for determining the off-diagonal elements of the MLF and the operator equations (23) (or (24)) that define the diagonal element corresponding to a possibly unstable subsystem. These results develop the results of [24] for a wider class of infinite-dimensional dynamical systems.

% Наведені в розділі 6 приклади показують, що запропонований метод побудови функцій Ляпунова застосовний у випадку, коли одна з підсистем є нестійкою, а отже відомі методи дослідження стійкості зв'язаних систем засновані на ВФЛ чи small-gain теоремах не застосовні. Також на прикладах показано, що для випадку коли обидві підсистеми експоненціально стійкі і обидва підходи застосовні, наш підхід призводить до менш консервативних умов стійкості. Значною перевагою запропонованого підходу є його простота, адже фактично задача про стійкість зв'язаної системи зводиться до розв'зання операторного рівняння, яке містить лише обмежені оператори і зводиться до матричного рівняння у випадку $\dim H_1<\infty$. Отже запропонований нами підхід значно розширює можливості методу функцій Ляпунова для широкого класу нескінченно-вимірних систем. Певним обмеженням нашого підходу є умова відділення спектрів \eqref{3.11}. В подальших дослідженнях планується послабити цю умову або відмовитися від неї. Можлива альтернатива умові \eqref{3.11} припущення про цілком неперервність лінійних операторів $A$ і $D^{-1}$, яке є природнім у випадку, коли друга підсистема є параболічною. У цьому випадку виникає можливість використання теорії Ріса---Шаудера цілком неперервних операторів в гільбертовому просторі. Інший напрямок узагальнення отриманих результатів це розгляд лінійних систем у банахових просторах, що дозволить досліджувати $L^p$- стійкість.
The examples given in section 6 show that the proposed method of constructing Lyapunov functions is applicable in the case when one of the subsystems is unstable, and therefore the known methods of studying the stability of coupled systems based on VLF or small-gain theorems are not applicable. The examples also show that for the case when both subsystems are exponentially stable and both approaches are applicable, our approach leads to less conservative stability conditions. A significant advantage of the proposed approach is its simplicity, since in fact the problem of the stability of a coupled system is reduced to solving an operator equation which contains only limited operators and is reduced to a matrix equation in the case $\dim H_1<\infty$. Therefore, our proposed approach significantly expands the capabilities of the Lyapunov function method for a wide class of infinite-dimensional systems. A certain limitation of our approach is the separation of spectrum condition \eqref{3.11}. In further research, it is planned to weaken this condition or abandon it. A possible alternative to the condition \eqref{3.11} is the assumption of complete continuity of the linear operators $A$ and $D^{-1}$ which is natural in the case when the second subsystem is parabolic. In this case, it becomes possible to use the Ries---Schauder theory of completely continuous operators in Hilbert space. Another direction of generalization of the obtained results is the consideration of linear systems in Banach spaces which allows to study $L^p$-stability.

\begin{ack}                               % Place acknowledgements
  % here.
\end{ack}

\bibliographystyle{plain}        % Include this if you use bibtex 
\bibliography{autosam}           % and a bib file to produce the 
                                 % bibliography (preferred). The
                                 % correct style is generated by
                                 % Elsevier at the time of printing.

\appendix

\section{Proof of Theorem 3} Assumption 3 implies Assumption 1 with setting $P_2 = \id_{H_2}$. Theorem 1 then provides the well-posedness of \eqref{3.46*}.\\ 
We show first the following lema:
\begin{lemma}
Let $H$ be a Hilbert space of dimension at least 2, $a,b\in H$. Then $J:=a\otimes b+b\otimes a$ is self-adjoint and 
\begin{equation*}
\gathered
\lambda_{\max}(J)
=
\Ree\langle a,b\rangle_{H}+\sqrt{\|a\|_H^2\|b\|_H^2-(\Imm\langle a,b\rangle_{H})^2}
\endgathered
\end{equation*}
\end{lemma}
{\it Proof of Lemma 3.} If $a$ and $b$ are linearly dependent, the statement can be proven by using the Cauchy-Buniakovsky inequality. Suppose that $a$ and $b$ are linearly independent, then $\{a\}^{\bot} \neq \{b\}^{\bot}$ and from
\begin{equation*}
\gathered
Jx
=\langle b,x\rangle_{H}a+\langle a,x\rangle_{H}b,
\endgathered
\end{equation*}
it follows $J \in L(H), J^{\dag}=J, R(J)=\spann\{a,b\}$. We decompose $x$ with respect to $H=R(J)\oplus R^{\bot}(J)$ and get $x=x_{\|}+x_{\bot}$, $x_{\|}\in R(J)$, $x_{\bot}\in R^{\bot}(J)$.
Then
\begin{equation}\label{Ineq1}
\gathered
\lambda_{\max}(J)=\sup\limits_{x\neq 0}
\frac{\langle x,Jx\rangle_{H}}{\|x\|_H^2} = \sup\limits_{x \neq 0}\frac{\langle x_{\|},Jx_{\|}\rangle_{H}}{\|x_{\bot}\|^2_H+\|x_{\|}\|^2_H}\\
=\max\Big\{0,\lambda_{\max}(J\big|_{R(J)})\Big\}.
\endgathered
\end{equation}
Rayleigh-Ritz inequality provides $\lambda_{\max}(J|_{R(J)})=\max_{\lambda\in\sigma(J|_{R(J)})}\lambda$.

Matrix of a linear operator $J|_{R(J)}$  with respect to a basis $\{a,b\}$ is 
\begin{equation*}
\gathered
\mathcal J=\begin{pmatrix}
\langle b,a\rangle_{H}&\|b\|_H^2\\
\|a\|_H^2&\langle a,b\rangle_{H}
\end{pmatrix}
\endgathered
\end{equation*}
By direct calculations
\if 0
\begin{equation*}
\gathered
\sigma(J|_{R(J)})=\sigma(\mathcal J)\\=
\Big\{\Ree\langle a,b\rangle_{H}\pm\sqrt{\|a\|_H^2\|b\|_H^2-(\Imm\langle a,b\rangle_{H})^2}\Big\}
\endgathered
\end{equation*}
\fi
\begin{equation*}
\gathered
\max_{\lambda \in \s (\mathcal J)}(\lambda)=\Ree\langle a,b\rangle_{H}+\sqrt{\|a\|_H^2\|b\|_H^2-(\Imm\langle a,b\rangle_{H})^2}\ge 0.
\endgathered
\end{equation*}

The statement of the lemma follows from $\sigma(J|_{R(J)})=\sigma(\mathcal J)$ and \eqref{Ineq1}.

{\it Proof of Theorem 3.}\\
For the proof of the theorem, we move according to the plan of section 4, defining the scalar Lyapunov function 
\begin{equation*}
    v(x,y)=p_{11}|x|^2+2\Ree(\overline{x}P_{12}y)+\|y\|^2
\end{equation*}
with $p_{11}\in \mathbb R_{>0}, P_{12} \in L(H,\mathbb C)$.\\
Analogically to \eqref{P_12_EQ}, we choose the linear operator $P_{12}$ as a solution of the operator equation
\begin{equation*}
\gathered
\overline a P_{12}+P_{12}D=-(p_{11}B
+C^{\dag}),
\endgathered
\end{equation*}
where $B\,:\,H \to \mathbb C$, $y\mapsto \langle b,y \rangle$ and $C:\mathbb C\to H$, $x\mapsto cx$.
Therefore, $P_{12}:H \to \mathbb C$ is defined by
\begin{equation}\label{P_12_EQ Appendix}
\gathered
P_{12} y = \langle R^{\dag}_D(- \overline a)(p_{11}b
+c),y \rangle 
\endgathered
\end{equation}
and we calculate the norm of it as follows:
\begin{equation*}
\gathered
\|P_{12}\|=\|R^{\dag}_D(- \overline a)(p_{11}b
+c)\|.
\endgathered
\end{equation*}
Thus,
\begin{equation*}
\gathered
\|P_{12}\|^2=p_{11}^2\rho_b^2+2p_{11}\Ree\mu_{bc}+\rho_c^2.
\endgathered
\end{equation*}
To find conditions for the coerciveness of the Lyapunov function we use the estimates
\begin{equation*}
\gathered
v(x,y)\geq p_{11}|x|^2-2\|P_{12}\||x|\|y\|+\|y\|^2\\\geq \lambda_{\min}(U_-)(|x|^2+\|y\|^2)\\
v(x,y)\leq p_{11}|x|^2+2\|P_{12}\||x|\|y\|+\|y\|^2\\\leq \lambda_{\max}(U_+)(|x|^2+\|y\|^2)
\endgathered
\end{equation*} with
\begin{equation*}
U_{\pm} =
\begin{pmatrix}
p_{11}&\pm\|P_{12}\|\\
\pm\|P_{12}\|&1
\end{pmatrix}.   
\end{equation*}
For coerciveness we need $\lambda_{\min}(U_-)=\min_{\lambda \in \s (U_-)}(\lambda)$ to be positive, which is equivalent to condition
\begin{equation*}
\gathered
\rho_b^2p_{11}^2+(2\Ree\mu_{bc}-1)p_{11}+\rho_c^2 < 0.
\endgathered
\end{equation*}
Using $2\Ree \mu_{bc}<1$ and $4\chi_b\chi_c<1$ we find the condition on $p_{11}$:
\begin{equation*}
\gathered
p_{11}\in(\theta_-,\theta_+).
\endgathered
\end{equation*}
Further, we investigate $\dot v(x,y)$ for $y \in \mathcal D(D)$. Due to the choice of $P_{12}$,
\begin{equation}\label{Dot_v_eq_1}
    \gathered
    \dot v(x,y)=\overline{x}\big((2p_{11}\Ree a)\id_{\mathbb C}+(P_{12}C+C^{\dag}P_{12}^{\dag})\big)x\\+2\Ree\langle y,Dy\rangle+\langle y,(B^{\dag}P_{12}+P_{12}^{\dag}B)y\rangle
    \endgathered
\end{equation} holds.
By direct calculations, from \eqref{P_12_EQ Appendix} we get 
\begin{equation*}
    P_{12}C+C^{\dag}P_{12}^{\dag}=2\Ree(p_{11}\rho_{bc}+\rho_{cc})\id_{\mathbb C}.
\end{equation*} 
and
\begin{equation*}
    \gathered
    B^{\dag}P_{12}+P_{12}^{\dag}B=p_{11}(b \otimes g_b+g_b \otimes b)\\+(b \otimes g_c + g_c \otimes b),
    \endgathered
\end{equation*} where
\begin{equation*}
g_b=R^{\dag}_D(-\overline{a})b, \,\, g_c=R^{\dag}_D(-\overline{a})c.
\end{equation*}
\if 0
Setting the $x$-term in \eqref{Dot_v_eq_1} as $-2\id_{\mathbb C}$, we imply
\begin{equation*}
    \gathered
    p_{11}=-\frac{1+\Ree\rho_{cc}}{\Ree(a+\rho_{bc})}.
    \endgathered
\end{equation*}
Substituting it in \eqref{Eta_Cond_2} we get \eqref{SC2}. Condition \eqref{SC1} is equivalent to $p_{11}$ being positive.\\ 
\fi
Due to Assumption 3 and Lemma 3, we estimate \eqref{Dot_v_eq_1} as follows:
\begin{equation*}
\gathered
\dot v(x,y) \leq 2\big(\Ree(a+\rho_{bc})p_{11}+\Ree\rho_{cc}\big)|x|^2\\+\big(-2\gamma_0+p_{11}\omega_{bb}+\omega_{bc}\big)\|y\|^2
\endgathered
\end{equation*}
\if 0
Condition \eqref{SC3} is equivalent to inequality
\begin{equation*}
(2\gamma_0-\omega_{bc})-p_{11}\omega_{bb}>0.
\end{equation*}
\fi
Sufficient conditions on $v$ to be a Lyapunov function are:
\begin{equation*}
\gathered
p_{11}\in(\theta_-,\theta_+), \quad \omega_{bb}p_{11}+\omega_{bc}<2\gm_0,\\
\Ree(a+\rho_{bc})p_{11}+\Ree\rho_{cc}<0.
\endgathered
\end{equation*}
Existence of parameter $p_{11}>0$ is assured by the statement of the theorem.
Exponential stability follows in the same way as in steps \eqref{3.32}-\eqref{3.33}.

Theorem 3 is proved.

\section{Proof of Proposition 1} 
Here we analyse the proof of Theorem 3, paying attention to the fact, that we can only obtain the estimates on parameters introduced there, since we cannot calculate the resolvent $R_D(-a)$ explicitly.

Operator $D^{-1}$ is self-adjoint and so is the resolvent $R_D(-a)$, the existence of which is assured by the domination condition \eqref{3.11}\,:
\begin{equation*}
\gathered
R_D^{\dag}(-a)=R_{D}(-a)=-D^{-1}(\id+a D^{-1})^{-1}\\=\sum\limits_{k=0}^{\infty}(-1)^{k+1}a^kD^{-k-1}
\endgathered
\end{equation*} with the norm-estimate
\begin{equation*}
\gathered
\|R_D(-a)\|\le\sum\limits_{m=0}^{\infty}\Big\|\sum\limits_{k=0}^{r-1}(-1)^{mr+k+1}a^{mr+k}D^{-mr-k-1}\Big\|\\
\le \Big\|\sum\limits_{k=0}^{r-1}(-1)^{k+1}a^{k}D^{k-1}\Big\|\sum\limits_{m=0}^{\infty}\big(|a|^r\|D^{-r}\|\big)^{m}\\= \frac{1}{1-q_r}\Big\|\sum\limits_{k=0}^{r-1}(-1)^{k+1}a^{k}D^{-k-1}\Big\|<\infty,
\endgathered
\end{equation*}
Analogically, we estimate all parameters introduced in Theorem 3 with upper-bounds $(\,{\cdot}\,)^{(r)}$.\\
\if 0
2){\bf Estimates of the introduced parameters:}\\2.1)
\begin{equation*}
\gathered
\rho_b=\|R_D^{\dag}(-a)b\|\le{\rho}_b^{(r)}
\endgathered
\end{equation*}
\begin{equation*}
\gathered
\rho_b^{(r)}=\Big(\int\limits_0^{\ell}\Big(\int\limits_0^{\ell}\Psi_r(z,\zeta)b(\zeta)\,d\zeta\Big)^2\,dz\Big)^{1/2}\\
\Psi_r(z,\zeta)=\sum\limits_{k=0}^{r-1}(-1)^{k+1}a^k\Gamma_{k+1}(z,\zeta),\\
%\delta_b^{(r)}:=\frac{|a|^{r+1}}{1-q}\Big(\int\limits_0^{\ell}\Big(\int\limits_0^{\ell}\Gamma_{r+2}(z,\zeta)b(\zeta)\,d\zeta\Big)^2\,dz\Big)^{1/2}
\endgathered
\end{equation*}
2.2)
\begin{equation*}
\gathered
\rho_c
\le \frac{1}{1-q_r}\rho_c^{(r)}
=:{\rho}_c^{(r)}
\endgathered
\end{equation*}
\begin{equation*}
\gathered
\rho_c^{(r)}=\Big(\int\limits_0^{\ell}\Big(\int\limits_0^{\ell}\Psi_r(z,\zeta)c(\zeta)\,d\zeta\Big)^2\,dz\Big)^{1/2}\\
%\delta_c^{(r)}:=\frac{|a|^{r+1}}{1-q}\Big(\int\limits_0^{\ell}\Big(\int\limits_0^{\ell}\Gamma_{r+2}(z,\zeta)c(\zeta)\,d\zeta\Big)^2\,dz\Big)^{1/2}
\endgathered
\end{equation*}
2.3)
\begin{equation*}
\gathered
\mu_{bc}\le\mu_{bc}^{(r)}+q_r(2-q_r) \rho_b^{(r)} \rho_c^{(r)}=:{\mu}_{bc}^{(r)},
\endgathered
\end{equation*}
\begin{equation*}
\gathered
\mu_{bc}^{(r)}=\int\limits_0^{\ell}\int\limits_0^{\ell}\int\limits_0^{\ell}\Psi_r(z,\zeta)
\Psi_r(z,\xi)b(\zeta)c(\xi)\,dz\,d\xi\,d\zeta.
\endgathered
\end{equation*}
{\bf Changes in the proof of Theorem 3:}\\
\fi
Sufficient condition for coerciveness of scalar Lyapunov function $v(x,y)$ reads\,:
\begin{equation*}
\gathered
(\rho_b^{(r)})^2p_{11}^2+(2\Ree\mu_{bc}^{(r)}-1)p_{11}+( \rho_c^{(r)})^2 < 0, 
\endgathered
\end{equation*}which is equivalent to 
\begin{equation*}
\gathered
p_{11}\in(\theta_-^{(r)},\theta_+^{(r)}).
\endgathered
\end{equation*}
Two other conditions on $p_{11}$ we rewrite in terms of bounds of correspondent parameters\,:
\begin{equation*}
\gathered
\omega_{bb}^{(r)}p_{11}<\frac{2\pi^2 d_0}{\ell^2}-\omega_{bc}^{(r)},\quad
(a+\rho_{bc}^{(r)})p_{11}+\rho_{cc}^{(r)}<0.    
\endgathered
\end{equation*}
\if 0
\begin{equation*}
\gathered
\rho_{bc}\le{\rho}_{bc}^{(r)},\\
%\rho_{bc}^{(r)}=\int\limits_0^{\ell}\int\limits_0^{\ell}\Psi_{r}(z,\zeta)c(z)b(\zeta)\,dz\,d\zeta.
\endgathered
\end{equation*}
\begin{equation*}
\gathered
\rho_{cc}\le{\rho}_{cc}^{(r)},\\
%\rho_{cc}^{(r)}=\int\limits_0^{\ell}\int\limits_0^{\ell}\Psi_{r}(z,\zeta)c(z)c(\zeta)\,dz\,d\zeta.
\endgathered
\end{equation*}
\begin{equation*}
\gathered
\rho_{bb}\le{\rho}_{bb}^{(r)},\\
%\rho_{bb}^{(r)}=\int\limits_0^{\ell}\int\limits_0^{\ell}\Psi_{r}(z,\zeta)b(z)b(\zeta)\,dz\,d\zeta.
\endgathered
\end{equation*}
\begin{equation*}
\gathered
\omega_{bb}=\rho_{bb}+\|b\|\rho_b\leq{\rho}_{bb}^{(r)}+\|b\|{\rho}_{b}^{(r)}={\omega}_{bb}^{(r)},\\
\omega_{bc}=\rho_{bc}+\|b\|\rho_c\leq{\rho}_{bc}^{(r)}+\|b\|{\rho}_{c}^{(r)}={\omega}_{bc}^{(r)}.
\endgathered
\end{equation*}
\fi
The received inequalities for parameter $p_{11}$ are equivalent to those in the proof of Theorem 3, so the investigation of solvability of them with respect to $p_{11}$ leads to the system of inequalities in Theorem 3, rewritten with upper-bounds $(\,{\cdot}\,)^{(r)}$.\\
%The last inequality is equivalent to the system of inequalities \eqref{3.66}, which completes the proof of Proposition 3.
\\
The Proposition is proved.
\end{document}